\documentclass[11pt]{article}
\usepackage{amsfonts}

\usepackage{graphics}
\usepackage{indentfirst}
\usepackage{cite}
\usepackage{latexsym}
\usepackage{amsmath}
\usepackage{amssymb}
\usepackage[dvips]{epsfig}
\usepackage{amscd}

\newtheorem{theorem}{Theorem}[section]
\newtheorem{remark}{Remark}[section]

\newtheorem{definition}{Definition}[section]
\newtheorem{lemma}[theorem]{Lemma}

\newtheorem{proposition}[theorem]{Proposition}

\newcommand{\n}{\rho}
\newcommand{\nr}{\tilde\rho^R}
\newcommand{\ti}{\tilde}

\newcommand{\lm}{\lambda}
\newcommand\divg{{\text{div}}}
\def\pf{{\it Proof.}  }
\renewcommand{\div}{ {\rm div }  }

\newcommand{\na}{\nabla }
\newcommand{\vp}{\varphi }

\newcommand{\pa}{\partial}

\newcommand{\bi}{\bibitem}

\newcommand{\bt}{\begin{theorem}}
\newcommand{\bl}{\begin{lemma}}
\newcommand{\el}{\end{lemma}}
\newcommand{\et}{\end{theorem}}
\newcommand{\ga}{\gamma}

\newcommand{\rot}{{\rm rot} }
\newcommand{\te}{\theta}

\newcommand{\al}{a }
\newcommand{\de}{\delta}
\newcommand{\ve}{\varepsilon}
\newcommand{\la}{\label}

\newcommand{\bn}{\begin{eqnarray}}
\newcommand{\en}{\end{eqnarray}}
\newcommand{\bnn}{\begin{eqnarray*}}
\newcommand{\enn}{\end{eqnarray*}}

\newcommand{\bnnn}{\begin{eqnarray*}}
\newcommand{\ennn}{\end{eqnarray*}}
\newcommand{\ben}{\begin{enumerate}}
\newcommand{\een}{\end{enumerate}}

\newcommand{\ba}{\begin{aligned}}
\newcommand{\ea}{\end{aligned}}
\newcommand{\be}{\begin{equation}}
\newcommand{\ee}{\end{equation}}

\def\p{\partial}
\def\norm[#1]#2{\|#2\|_{#1}}

\def\lap{\triangle}

\def\lam{\lambda}
\def\ep{\varepsilon}

\def\o{\omega}

\def\rr{\mathbb{R}^2}
\def\O{B_R}

\makeatletter      
\@addtoreset{equation}{section}
\makeatother       

\title{On    Classical Solutions   to the Cauchy Problem of the Two-Dimensional Barotropic
Compressible Navier-Stokes Equations with Vacuum}
\date{}
\author{Jing L{\small I},\thanks{Institute of Applied Mathematics, AMSS, and Hua Loo-Keng Key Laboratory of Mathematics,
Chinese Academy of Sciences, Beijing 100190, People's  Republic of
China ({\tt ajingli@gmail.com}).
 The research of \textsc{J. Li} was
partially supported  by the National Center for Mathematics and Interdisciplinary Sciences, CAS, and NNSFC   11171326.}  \quad Zhilei L{\small IANG} \thanks{School of Economic Mathematics, South Western University of Finance  and Economics, Chengdu 611130, People's  Republic of
China  ({\tt zhilei0592@gmail.com}). The research of \textsc{Z. Liang} was
partially supported  by  NNSFC   11226163. }
  }

\begin{document}
\maketitle

\begin{abstract}This paper concerns  the
Cauchy problem of the barotropic compressible Navier-Stokes equations  on the whole two-dimensional space  with vacuum as far field density. In particular, the initial density can have compact support.
When the shear and the bulk  viscosities are  a positive constant and  a
power function of the density respectively,  it is proved that   the  two-dimensional  Cauchy problem of
the compressible Navier-Stokes equations   admits a unique local strong
solution provided the initial density decays not too slow at infinity. Moreover, if the initial data satisfy  some additional regularity  and   compatibility   conditions, the strong solution becomes a classical one.
\end{abstract}


\section{Introduction and main results}
  We consider the two-dimensional
  barotropic compressible Navier-Stokes equations    which read as follows:
\be\la{n1}
\begin{cases} \rho_t + \div(\rho u) = 0,\\
 (\rho u)_t + \div(\rho u\otimes u) + \nabla P = \mu\lap u + \nabla((\mu + \lam)\div u),
\end{cases}
\ee
where   $t\ge 0, x=(x_1,x_2)\in  \Omega\subset \rr, \rho=\n(x,t)$ and $u=(u_1(x,t),u_2(x,t))$ represent, respectively,  the density and the velocity,  and the pressure $P$ is given by
\be\la{n2}
P(\rho) = A\rho^{\gamma}\, ( A>0,\, \ga>1).
\ee
The shear viscosity   $\mu$ and the bulk one $\lambda$ satisfy  the following hypothesis:
\be\la{n3} 0<
\mu = const,\quad \lam(\n)  = b\rho^{\beta},
\ee where the constants $b$ and $\beta\ge 0$ satisfy \be \begin{cases}b>0,&\mbox{ if } \,\beta>0,\\ \mu +b\ge 0,&\mbox{ if }\, \beta=0. \end{cases} \ee
 In the sequel, without loss of generality, we set $A= 1.$
 Let $\Omega=\rr$ and we consider the Cauchy problem with $(\n,u)$ vanishing at infinity (in some weak sense). For given initial data $\n_0$ and $m_0,$   we require that
\be \la{n4} \n(x,0)=\n_0(x), \quad \n u(x,0)=m_0(x),\quad x\in \Omega=\rr.\ee

When both the shear and bulk viscosities are positive constants, there are extensive studies  concerning the theory of strong and weak solutions for the system of the multi-dimensional compressible Navier-Stokes equations.
When the data $\n_0,m_0$ are sufficiently regular and
the initial density $\n_0$ has a positive lower bound, there exist  local strong and  classical solutions
to the problem \eqref{n1} and the solutions exist  globally in time provided that the data
are small in some sense. For details, we refer the readers to  \cite{Na,se1,M1,dan,Ho4} and the references therein. On the other hand, in the case that the initial density need not be positive
and may vanish in open sets, under some additional compatibility conditions, the authors in \cite{cho1,
K2,coi1, sal}  obtained
the local existence and uniqueness of strong and
classical solutions   for  three-dimensional bounded or unbounded domains and for two-dimensional bounded ones. Later, the compatibility conditions on the initial data were further relaxed by Huang-Li-Matsumura \cite{hlma}. For  large initial data, the global existence of weak    solutions was first obtained by Lions \cite{L1,L4} provided $\ga$ is suitably large  which was  improved later by Feireisl-Novotny-Petzeltova \cite{F1}(see also \cite{Fe}). Recently, for the case that the initial density is allowed to
vanish,  Huang-Li-Xin \cite{hlx1} obtained the global existence  of classical
solutions to the Cauchy problem for the   isentropic compressible
Navier-Stokes equations in three spatial dimensions with smooth
initial data provided that the initial energy is suitably small. 

For  large  initial data away from vacuum, Vaigant-Kazhikhov \cite{Ka}   first obtained   that   the two-dimensional system  \eqref{n1}-\eqref{n4} admits a unique global strong solution   provided that $\beta>3$ and that $\Omega$ is bounded. Recently, for  periodic   initial data with initial density allowed to vanish,   Huang-Li\cite{hlia} relaxed the crucial condition $\beta>3$ of  \cite{Ka} to the one that  $\beta>4/3.$ However,  for the Cauchy problem  \eqref{n1}-\eqref{n4} with $\Omega=\rr,$ it is still open even for the local existence of strong and classical solutions when the far field   density is vacuum, in particular,  the initial density may have compact support. In fact, this is the aim of this paper.

In this section, for $1\le r\le \infty ,$  we  denote the standard Lebesgue and
  Sobolev spaces as follows:
$$ L^r=L^r(\rr), \quad W^{s,r}= W^{s,r}(\rr),  \quad H^s= W^{s,2} .  $$ The first main result of this paper  is  the  following Theorem \ref{t1} concerning the  local existence   of  strong solutions whose definition is as follows:
\begin{definition} If  all derivatives involved in \eqref{n1} for $(\rho,u)  $  are regular distributions, and   equations  \eqref{n1} hold   almost everywhere   in $\rr\times (0,T),$ then $(\n,u)$  is called a  strong solution to  \eqref{n1}.
\end{definition}

\begin{theorem}\la{t1} Let $\eta_0 $ be a positive constant and
  \be\la{2.07} \bar x\triangleq(e+|x|^2)^{1/2}\log^{1+\eta_0} (e+|x|^2) .\ee  For  constants $q>2 $ and $a\in (1,2) ,$ assume that  the initial data $( \n_0 ,m_0)$ satisfy that
  \be\la{1.9}\n_0\ge 0,\,\,
  \bar x^a \rho_0\in   L^1 \cap H^1\cap W^{1,q},\,\,
\na u_0 \in  L^2 , \,\,\rho_0^{1/2}u_0\in L^2,
  \ee and that\be \la{a1.9} m_0=\n_0u_0.\ee  In addition, if $\beta\in (0,1),$ suppose that \be\la{c1.9}\lm(\rho_0) \in L^2,\quad \bar x^{\theta_{0}} \na \lm(\rho_0) \in L^2\cap  L^{q},\ee for some $\theta_0\in(0,\min\{ \beta,1\}) .$ Then there exists a positive time $T_0>0$ such that the problem  \eqref{n1}-\eqref{n4} has a unique   strong solution $(\n,u)$ on $\rr\times (0,T_0]$ satisfying that\be\la{1.10}\begin{cases}
  \rho\in C([0,T_0];L^1 \cap H^1\cap W^{1,q} ),\quad  \bar x^a\rho\in L^\infty( 0,T_0 ;L^1\cap H^1\cap W^{1,q} ),\\ \sqrt{\n } u,\,\na u,\, \bar x^{-1}u,
   \,    \sqrt{t} \sqrt{\n}  u_t \in L^\infty(0,T_0; L^2 ) , \\ \na u\in  L^2(0,T_0;H^1)\cap  L^{(q+1)/q}(0,T_0; W^{1,q}), \,\, \sqrt{t}\na u\in L^2(0,T_0; W^{1,q} )   ,  \\ \sqrt{\n} u_t, \, \sqrt{t}\na u_t ,\,  \sqrt{t} \bar x^{-1}u_t\in L^2(\rr\times(0,T_0)) , \\ \lm(\n)\in C([0,T_0];L^2 ),\,\, \bar x^{\te_0}\na \lm(\n)  \in L^\infty( 0,T_0 ;L^2\cap L^q ) ,\,\,\mbox{ for }\beta\in (0,1),
   \end{cases}\ee   and that \be \la{l1.2}\inf\limits_{0\le t\le T_0}\int_{B_{N }}\n(x,t) dx\ge \frac14\int_{\rr} \n_0(x)dx,\ee
for some constant $N >0$ and $B_{N }\triangleq\left.\left\{x\in\rr\right|
\,|x|<N \right\} .$
\end{theorem}

Furthermore, if the initial data $(\n_0,m_0)$ satisfy some additional regularity and compatibility conditions, the local strong solution $(\n,u) $ obtained by Theorem \ref{t1} becomes a classical one for positive time, that is, we have
\begin{theorem}\la{t2} In addition to    \eqref{1.9}-\eqref{c1.9},  assume further that
  \be\la{1.c1}
  \begin{cases}  \na^2 \n_0,\,\na^2 \lm(\n_0 ),\,\na^2 P(\n_0 )\in L^2\cap  L^q,\\ \bar x^{\de_0}\na^2   \rho_0 ,\, \bar x^{\de_0}\na^2  \lm(\rho_0  ),\,\bar x^{\de_0} \na^2 P(  \rho_0)   \in  L^2    , \quad  \na^2 u_0 \in  L^2 ,
   \end{cases}\ee   for some     constant $\de_0\in(0,\te_0]. $ Moreover, suppose  that  the following  compatibility conditions hold   for  some   $g\in L^2,  $
\be \la{co2}- \mu\lap u_0 - \nabla((\mu + \lam(\n_0))\div u_0)+  \nabla P(\n_0)=\n_0^{1/2}g . \ee      Then,   in addition to \eqref{1.10} and  \eqref{l1.2},    the  strong  solution  $(\rho,u)$ obtained by Theorem \ref{t1} satisfies   \be\la{1.a10}\begin{cases}
  \na^2\rho, \,\, \na^2\lm(\rho), \,\,\na^2 P(\rho)\in C([0,T_0];L^2\cap L^q  ), \\ \bar x^{\de_0}\na^2  \rho   ,\,\, \bar x^{\de_0} \na^2  \lm(\rho )   ,\,\, \bar x^{\de_0} \na^2  P(  \rho)     \in L^\infty( 0,T_0 ;L^2 ) ,\\     \na^2 u,\, \sqrt{\n}  u_t, \,\sqrt{t}   \na  u_t,\,\sqrt{t}  \bar x^{-1}  u_t,\, t\sqrt{\n}u_{tt}, \,t  \na^2 u_t\in L^\infty(0,T_0; L^2),\\ t\na^3 u\in  L^\infty(0,T_0; L^2\cap L^q), \,\\    \na u_t,\, \bar x^{-1}u_t,\,  t\na u_{tt},\, t\bar x^{-1}u_{tt}\in L^2(0,T_0;L^2), \\   t \na^2(\n u)\in L^\infty(0,T_0;L^{(q+2)/2}) .
   \end{cases}\ee    \end{theorem}

A few remarks are in order:

\begin{remark} First, it follows from (\ref{1.10}) and  (\ref{1.a10}) that
\be \la{1.z11} \n ,\,\, \bar x\n ,\,\, \na\n \in C(\overline{\rr} \times[0,T_0] ),\ee
and that \be\la{1.z1}\ba \n u\in  H^1(0,T_0; L^2)\hookrightarrow  C( [0,T_0]; L^2 ).\ea\ee
 The  Gargilardo-Nirenberg  inequality   shows that for $k=0,1,$
\bnn \|\na^k(\n u)\|_{C(\overline{\rr})}\le C\|\n u\|_{L^2}^{((2-k)q-2k)/(3q+2)}\|\na^2(\n u)\|_{L^{(q+2)/2}}^{(k+1)(q+2)/(3q+2)}, \enn
which together with \eqref{1.z1} and  \eqref{1.a10}  yields that for any  $ \tau\in (0,T_0 ),$
\be\la{1.z8} \n u  \in  C( [\tau,T_0];C^1( \overline{\rr}) ),\,\,\n_t  \in  C(\overline{\rr} \times[\tau,T_0] ).\ee

Next, we deduce from  (\ref{1.10}) and  (\ref{1.a10}) that
\be  \bar x^{-1}  u  \in L^\infty( 0,T_0; H^2)\cap H^1(0,T_0;L^2) \hookrightarrow  C(\overline{\rr} \times[0,T_0]),\ee and that for any  $ \tau\in (0,T_0 ),$
\be  \la{1.z9}\na u  \in H^1(0,T_0;L^2)\cap L^\infty( 0,T_0; H^1) \cap L^\infty(\tau,T_0;  W^{2,q})\hookrightarrow  C([\tau,T_0] ;C^1(\overline{\rr} ) ).\ee

Finally, for any  $ \tau\in (0,T_0 ),$ it follows  from  (\ref{1.10}) and  (\ref{1.a10})  that
\be  \la{1.z10} \na u_t,\,\, \bar x^{-1}u_t  \in   H^1(\tau,T_0;L^2) \hookrightarrow  C([\tau,T_0] ;L^2 ),\ee
which combined with $\na^2 u_t \in L^\infty(\tau,T_0;L^2)$ gives
 \be  \la{1.zz10} \na u_t  \in   C([\tau,T_0] ;L^p ),\quad \bar x^{-1}u_t\in  C(\overline{\rr} \times[\tau,T_0] ), \ee for any $p\ge 2.$ This together with \eqref{1.z11} shows \bnn \n u_t\in  C(\overline{\rr} \times[\tau,T_0] ),\enn
    which, along with \eqref{1.z11} and \eqref{1.z8}-\eqref{1.zz10},  thus implies that  the
   solution $(\n,u )$ obtained by Theorem \ref{t1} is  in fact a classical
    one to  the Cauchy problem  (\ref{n1})-(\ref{n4})
     on $\rr\times (0,T_0].$
 \end{remark}

\begin{remark}  To obtain the local existence and uniqueness of strong solutions, in Theorem \ref{t1}, the only compatibility condition we   need is \eqref{a1.9} which is  similar to that of \cite{hlma} and is much weaker than those of   \cite{cho1,
K2,coi1, sal}  where   not only \eqref{a1.9} but also \eqref{co2}  is   needed. Moreover, for the local existence   of classical solutions, Cho-Kim \cite{coi1}  needs the following additional  condition:
  \bnn \na (\n_0^{-1/2}g)\in L^2(\rr), \enn besides   \eqref{a1.9} and \eqref{co2}. This is in fact stronger than the compatibility conditions listed in our Theorem \ref{t2}.
 \end{remark}

We now comment  on the analysis of this paper.   As mentioned by \cite{cho1,K2,coi1,hlma}, the methods in  \cite{cho1,
K2,coi1,hlma}  can not be applied directly to  our case, since  for two-dimensional case   their arguments   only works for the case that $\Omega$ is bounded.  In fact, for  $\Omega=\rr,$ it seems   difficult to bound the $L^p(\rr)$-norm of $u$ just in terms of $\|\n^{1/2} u\|_{L^2(\rr)} $ and   $\|\na u\|_{L^2(\rr)}.$   The key observations    to overcome the difficulties caused by the unbounded domain  are as follows:  On the one hand,  for system  \eqref{n1}, it is enough to bound the   $L^p(\rr)$-norm of the momentum $\n u$ instead of just the velocity $u,$  and on the other hand, since $\n$ decays for large values of the spatial variable $x,$  the momentum $\n u$   decays faster than $u$ itself.
To this end, we first establish a key Hardy-type inequality (see \eqref{3.a2})  by combining a Hardy-type one due to Lions \cite{L2} (see \eqref{2.i2}) with a spatial weighted mean estimate  of the density (see \eqref{p1}).
  We then construct the approximate solutions to  \eqref{n1}, that is, for density strictly away from vacuum initially, we consider \eqref{n1} in any bounded ball $B_R$   with radius $R>0.$ To overcome the difficulties caused by the fact that the bulk viscosity $\lm$  depends on $\n,$ we imposed the Navier-slip boundary conditions on \eqref{n1} instead of the usual Dirichlet boundary ones. However,  when we extend the approximate solutions by $0$ outside the ball, it seems difficult to bound the $L^2(\rr)$-norm of the gradient of the velocity. This will  be  overcome by  putting an additional term $-R^{-1}u$ on the right-hand side of $\eqref{n1}_2.$  See \eqref{b2} for details.  Finally, combining  all these ideas stated above with those due to  \cite{cho1,
K2,coi1,hlma}, we derive some desired   bounds on the gradients of the velocity and the   spatial weighted  ones on both the density and its gradients where all these bounds are independent of both the radius of the balls $B_R$ and the lower bound  of
the initial density.

The rest of the paper is organized as follows: In Section 2, we collect some
elementary facts and inequalities which will be needed in later analysis. Sections 3  and 4
are devoted to the a priori estimates   which  are  needed to obtain  the local existence and uniqueness of stong and classical solutions. Then finally, the main results,
Theorems \ref{t1} and \ref{t2}, are proved in Section 5.

\section{Preliminaries}

First, the following    local existence theory on bounded balls, where the initial
density is strictly away from vacuum, can be shown by similar arguments as in   \cite{cho1,K2,coi1}.
\begin{lemma}   \la{th0} For $R>0$  and $B_R=\{x\in\rr ||x|<R\},$ assume  that
 $(\n_0,u_0 )$ satisfies \be \la{2.1}
 (\n_0,u_0) \in H^3(\O), \quad    \inf\limits_{x\in\O}\n_0(x) >0 , \quad  u_0\cdot n=0,
 \,\,\rot u_0=0,\, x\in\pa\O.\ee   Then there exist  a small time
$T_R>0$ and a unique classical solution $(\rho , u )$ to the following initial-boundary-value problem   \be \la{b2} \begin{cases}\rho_t + \div(\rho u) = 0,\\
 (\rho u)_t + \div(\rho u\otimes u) + \nabla P - \mu\lap u - \nabla((\mu + \lam)\div u)=-R^{-1}u,   \\ u\cdot n =0,\, \rot u=0,\quad   x\in\pa\O, \,\, t>0,\\ (\n,   u)(x,0)=(\n_0,u_0)(x),\quad   x\in \O,\end{cases}\ee
  on
$\O\times(0,T_R]$ such that
\be\la{lo1}\begin{cases}    \n\in C\left([0,T_{R}];H^{3}\right), u\in C\left([0,T_{R}]; H^{3}\right)\cap L^{2}\left(0,T_{R};H^{4}\right),\\
 u_{t}\in L^{\infty}\left(0,T_{R};H^{1}\right)\cap  L^{2}\left(0,T_{R};H^{2}\right),\sqrt{\n}u_{tt}\in L^{2}\left(0,T_{R};L^{2}\right),\\
 \sqrt{t}  u \in L^{\infty}\left(0,T_{R};H^{4}\right), \, \sqrt{t} u_{t}\in  L^{\infty}\left(0,T_{R}; H^{2}\right),\,\sqrt{t} u_{tt}\in L^{2}\left(0,T_{R};H^{1}\right),\\
 \sqrt{t}\sqrt{\n}u_{tt}\in L^{\infty}\left(0,T_{R};L^{2}\right), tu_{t}\in L^{\infty}\left(0,T_{R};H^{3}\right), \\
  tu_{tt}\in L^{\infty}\left(0,T_{R};H^{1}\right) \cap  L^{2}\left(0,T_{R};H^{2}\right), t\sqrt{\n}u_{ttt}\in L^{2}\left(0,T_{R};L^{2}\right),\\
 t^{3/2} u_{tt}\in  L^{\infty}\left(0,T_{R};H^{2}\right),t^{3/2} u_{ttt}\in L^{2}\left(0,T_{R};H^{1}\right), \\ t^{3/2}\sqrt{\n}u_{ttt}\in L^{\infty}\left(0,T_{R};L^{2}\right),
 \end{cases}\ee where we denote $L^2=L^2(\O)$ and $H^k=H^k(\O)$ for positive integer $k .$
 \end{lemma}

   Next, for either $\Omega=\rr$ or $\Omega=B_R$ with $R\ge 1,$    the following weighted $L^p$-bounds for elements of the Hilbert space    $ \ti D^{1,2}(\Omega)\triangleq\{v\in H^1_{\rm loc}(\Omega)|\na v\in L^2(\Omega)\} $  can be found in \cite[Theorem B.1]{L2}.
\begin{lemma} \la{1leo}
   For   $m\in [2,\infty)$ and $\theta\in (1+m/2,\infty),$ there exists a positive constant $C$ such that   for  either $\Omega=\rr$ or $\Omega=B_R$ with $R\ge 1$ and for any  $v\in \ti D^{1,2}(\Omega) ,$ \be\la{3h} \left(\int_{\Omega} \frac{|v|^m}{e+|x|^2}(\log (e+|x|^2))^{-\theta}dx  \right)^{1/m}\le C\|v\|_{L^2(B_1)}+C\|\na v\|_{L^2 (\Omega)}.\ee
\end{lemma}

A useful consequence of Lemma \ref{1leo}  is the following weighted  bounds for elements of  $\ti D^{1,2}(\Omega) $ which in fact will play  a  crucial role   in our analysis.

\begin{lemma}\la{lemma2.6} Let  $\bar x$ and $\eta_0$ be as in \eqref{2.07} and $\Omega$ as in Lemma \ref{1leo}. For  $\ga>1,$ assume that $\n \in L^1(\Omega)\cap L^\ga(\Omega)$ is a non-negative function such that
\be \la{2.i2}   \int_{B_{N_1} }\n dx\ge M_1,  \quad \int_{\Omega} \n^\ga dx\le M_2,  \ee
for positive constants $   M_1, M_2, $ and $ N_1\ge 1$  with $B_{N_1}\subset\Omega.$ Then there is a positive constant $C$ depending only on   $   M_1,M_2,
 N_1, \ga,$ and $\eta_0$ such that     \be\la{3.ii2}\ba \|v\bar x^{-1}\|_{L^{2}(\Omega)} &\le C \|\n^{1/2}v\|_{L^2(\Omega)}+C \|\na v\|_{L^2(\Omega)} ,\ea\ee for   $v\in \ti D^{1,2}(\Omega)  .$  Moreover, for $\ve> 0$ and $\eta>0,$ there is a positive constant $C$ depending only on   $\ve,\eta, M_1,M_2,
 N_1, \ga,$ and $\eta_0$ such that      every $v\in \ti D^{1,2}(\Omega)  $ satisfies\be\la{3.i2}\ba \|v\bar x^{-\eta}\|_{L^{(2+\ve)/\ti\eta}(\Omega)} &\le C \|\n^{1/2}v\|_{L^2(\Omega)}+C \|\na v\|_{L^2(\Omega)} ,\ea\ee with $
\ti\eta=\min\{1,\eta\}.$

\end{lemma}

{\it Proof.}   It follows from   \eqref{2.i2} and  the Poincar\'e-type inequality     \cite[Lemma 3.2]{Fe}  that  there exists a positive constant $C $ depending only on $  M_1, M_2, N_1 ,$ and $\ga,$  such that  \bnn   \|v\|_{L^2(B_{ N_1} )}^2\le C \int_{B_{ N_1} }\n v^2dx +C \|\na v\|_{L^2(B_{ N_1} )}^2,\enn   which together with \eqref{3h} gives  \eqref{3.ii2} and \eqref{3.i2}.   The proof of Lemma \ref{lemma2.6} is finished.

Finally, the following $L^p$-bound for elliptic systems, whose proof is similar to that of \cite[Lemma 12]{cho1}, is a direct consequence of the combination of  a well-known elliptic theory due to Agmon-Douglis-Nirenberg \cite{adn} with a standard scaling procedure.
\begin{lemma} For $p>1$ and $k\ge 0,$ there exists a positive constant $C$ depending only on $p$ and $k$ such that
\be \la{lp} \|\na^{k+2}v\|_{L^p(\O)}\le C\| \Delta v\|_{W^{k,p}(\O)},\ee
for every    $v\in W^{k+2,p}(\O)$ satisfying either
\bnn v\cdot n=0,\,\, \rot v=0, \,\,\mbox{ on }\pa\O,\enn or \bnn v =0, \,\,\mbox{ on }\pa\O.\enn

\end{lemma}

\section{A priori estimates (I)}

Throughout this section and the next, for $p\in [1,\infty]$ and $k\ge 0,$
we    denote
\bnn \int fdx=\int_{\O}fdx, \quad L^p=L^p(\O),\quad W^{k,p}=W^{k,p}(\O),\quad H^k=W^{k,2},\enn and,
 without loss of generality, we assume that  $\beta>0 $ since all these estimates obtained in this section and the next hold for the case that $\beta=0$ after some small modifications.
Moreover, for   $R>4N_0\ge 4,$  assume that  $(\n_0,u_0)$ satisfies,  in addition to  \eqref{2.1}, that\be \la{w1} 1/2\le \int_{B_{N_0}}\n_0(x)dx \le \int_{B_R }\n_0(x)dx \le 3/2.\ee
  Lemma \ref{th0} thus yields that there exists some $T_R>0$ such that  the  initial-boundary-value problem   \eqref{b2} has a unique classical solution  $(\n,u)$ on $B_R\times[0,T_R]$ satisfying \eqref{lo1}.

 In this section, for  $\bar{x} $ and $\eta_0 >0$  as  in     \eqref{2.07} and  for $a\in(1,2),$ $q\in (2,\infty),$ and $\te_0>0$ as in Theorem \ref{t1}, we will use the convention that $C$ denotes a
generic positive constant
 depending only on       $\mu,\beta,\ga, b, q,  a, \eta_0, \te_0,  N_0, $ and $E_0,$ where  \bnn\ba E_0\triangleq &    \|\n_0^{1/2}u_0 \|_{L^2} + \|\na u_0\|_{L^2} +R^{-1/2}\|u_0\|_{L^2}+ \|\bar{x}^{a}\n_0\|_{L^{1}\cap H^{1}\cap W^{1,q}}\\&+\|\lm(\n_0)\|_{L^2}+\| \bar{x}^{\te_{0}}\na \lam(\n_0)\|_{L^2\cap L^q},\ea \enn and   we write $C(\kappa)$ to emphasize
that $C$ depends on $\kappa.$

 Denoting $ \na^{\bot}\triangleq(\p_{2},-\p_{1}),$ we rewrite the momentum equations $\eqref{b2}_{2}$ as
\be\la{pq1}\ba  \rho \dot{u} +R^{-1}u =  \na F+\mu\na^{\bot}  \o,\ea\ee where   $$ \dot f\triangleq
f_t+u\cdot\nabla f,\quad F\triangleq(2\mu + \lambda)\text{div}u -
P(\rho)  ,\quad\o \triangleq\na^\perp\cdot u, $$ are  the
material derivative of $f,$ the effective viscous flux, and the
vorticity respectively.
 Thus, \eqref{pq1} implies that $\o$ satisfies \be\la{a39} \begin{cases}  \mu\lap  \o =\na^{\bot}\cdot \left( \rho \dot{u}+R^{-1}u  \right),& \mbox{ in } B_R,\\ \o =0,& \mbox{ on }\pa B_R.\end{cases}\ee
Applying the standard $L^{p}$-estimate to \eqref{a39} yields that for $p\in (1,\infty),$
\bnn  \begin{cases} \| \na\o\|_{L^p}\le C(p)\left(\|\n \dot{u}\|_{L^{p}}  +R^{-1}\|u\|_{L^{p}} \right),\\ \| \na^2\o\|_{L^p}\le C(p)\left(\|\na(\n \dot{u})\|_{L^{p}}  +R^{-1}\|\na u\|_{L^{p}} \right) ,\end{cases}\enn
which together with \eqref{pq1} gives
\be   \la{2.9-1} \begin{cases}\| \na F \|_{L^p}+\| \na\o \|_{L^p}\le C(p)\left(\|\n \dot{u}\|_{L^{p}}+ R^{-1}\|u\|_{L^{p}} \right), \\  \| \na^2F \|_{L^p}+ \| \na^2\o\|_{L^p}\le C(p)\left(\|\na(\n \dot{u})\|_{L^{p}}  +R^{-1}\|\na u\|_{L^{p}} \right). \end{cases} \ee

The main aim of this section is to derive
  the following key a priori estimate  on $\psi $ defined by \be\ba\la{3.2}
\psi(t)\triangleq& 1 +\|\n^{1/2}u\|_{L^2} + \|\na u\|_{L^2} +R^{-1/2} \|  u\|_{L^2} \\&  + \|\bar{x}^{a}\n\|_{L^{1}\cap H^{1}\cap W^{1,q}}+\|\lm(\n)\|_{L^2}+\| \bar{x}^{\te_{0}}\na \lam(\n)\|_{L^2\cap L^q} .
\ea\ee
\begin{proposition} \la{pro}  Assume that $(\n_0,u_0)$ satisfies \eqref{2.1} and \eqref{w1}. Let $(\n,u)$ be the solution to  the  initial-boundary-value problem  \eqref{b2} on $\O\times (0,T_R]$ obtained by Lemma \ref{th0}. Then there exist positive constants $T_0$ and $M$ both  depending only on      $\mu,\beta,\ga,b, q, $ $  a, \eta_0, \te_0, N_0,$ and $E_0$ such that  \be\la{o1}\ba &\sup\limits_{0\le t\le T_0}\psi(t) +  \int_0^{T_0} \left(\norm[L^{q}]{\nabla^2u}^{(q+1)/q}+t\|\na^2u\|_{L^q}^{2}+ \|\na^2u\|_{L^2}^{2}\right) dt\le M.\ea\ee \end{proposition}

To prove Proposition \ref{pro} whose proof will be postponed to the end of this section, we begin with the following   standard   energy estimate for
$(\n,u)$ and preliminary    $L^2$-bounds for $\nabla u.$

\begin{lemma} \la{l3.0}  Let  $(\n,u)$ be a smooth solution to  the  initial-boundary-value problem  \eqref{b2}.  Then there exist a $T_1=T_1(E_0)>0$     and  a positive constant $\alpha=\alpha(\ga,\beta, q)>1$  such that for all $t\in (0, T_1],$
\be\ba\la{3.1}
  &\sup_{0\le s\le t} \left(R^{-1}\|u\|_{L^2}^2+\|\na u\|^2_{L^{2}}\right)
  +\int_{0}^{t}\int{\rho} |u_t|^2 dx ds   \le C  +C\int_{0}^{t}  \psi^{\alpha} ds.
\ea
\ee

\end{lemma}

\pf    First,  applying standard energy estimate to \eqref{b2} gives
\be\ba\la{3.1-c}
\sup\limits_{0\le s\le t}\left( \|\sqrt{\n}u\|_{L^{2}}^{2}+\|\n\|^\ga_{L^\ga}\right) +  \int_0^t\left(\|\na u\|_{L^{2}}^2  +R^{-1} \|  u\|_{L^{2}}^2\right)  ds
 \le C .
\ea\ee

Next,
for $N>1 $ and $\vp_N\in C^\infty_0(B_N)$   such that
 \be \la{vp1}0\le \vp_N \le 1, \quad  \vp_N(x)= 1, \mbox{ if } |x|\le N/2, \quad |\na^k \vp_N|\le C N^{-k} (k=1,2) ,\ee
  it follows from \eqref{3.1-c}  and \eqref{w1} that \be\la{oo0}\ba \frac{d}{dt}\int \n \vp_{2N_0} dx &=\int \n u \cdot\na \vp_{2N_0} dx\\ &\ge -  C N_0^{-1}\left(\int\n dx\right)^{1/2}\left(\int\n |u|^2dx\right)^{1/2}\ge - \ti C (E_0),\ea\ee  where in the last inequality we have used $$ \int\n dx =\int \n_0dx ,$$  due to $(\ref{b2})_1 $  and $(\ref{b2})_3 .$
  Integrating \eqref{oo0} gives  \be\la{p1}\ba \inf\limits_{0\le t\le T_1}\int_{B_{2N_0}}  \n dx&\ge \inf\limits_{0\le t\le T_1}\int \n \vp_{2N_0} dx\\ &\ge \int \n_0 \vp_{2N_0} dx-\ti C T_1  \ge 1/4,\ea\ee where $T_1\triangleq\min\{1, (4\ti C)^{-1}\} .$ From now on, we will always assume that $t\le T_1.$
  The combination of   \eqref{p1},     \eqref{3.1-c},    and  \eqref{3.i2}    yields that for $\ve> 0$ and $\eta>0 ,$  every $v\in   \ti D^{1,2}(\O) $  satisfies\be \la{3.v2}\ba  \|v\bar x^{-\eta}\|_{L^{(2+\ve)/\ti\eta} }^2 &\le C(\ve,\eta) \int \n |v|^2dx +C (\ve,\eta)  \|\na v\|_{L^2 }^2 ,\ea\ee with $\ti\eta=\min\{1,\eta\}.$
  In particular, we have
  \be\la{3.a2}\ba\|\n^\eta u\|_{L^{(2+\ve)/\ti\eta}}+ \|u\bar x^{-\eta}\|_{L^{(2+\ve)/\ti\eta}}  \le C(\ve,\eta)\psi^{1+\eta} .\ea\ee

Next, multiplying    equations  $\eqref{b2}_{2}$  by  $ u_t$ and    integration by parts yield
\be\la{3r1}\ba& \frac{d}{dt}\int\left((2\mu+\lm)(\div u)^2+\mu\o^2+R^{-1}|u|^2\right)dx+\int \n|  u_t|^2dx\\ &\le C\int\n|u|^2|\na u|^2dx+ \int \lm_t (\div u)^2dx+2\int P\div u_tdx .\ea\ee

We estimate each term on the right-hand side of \eqref{3r1} as follows:

First,   the Gagliardo-Nirenberg inequality implies  that for all $p\in(2,+\infty),$
\be \la{cc6}\ba\|\na u\|_{L^p}&\le C(p)\|\na u\|_{L^2}^{2/p}\|\na u\|_{H^1}^{1-2/p}\\&\le C(p)\psi +C(p)\psi\|\na^2 u\|_{L^2}^{1-2/p} ,\ea\ee which together with  \eqref{3.a2} yields that for $\eta>0 $ and $
\ti\eta=\min\{1,\eta\},$
\be\ba\la{cc5}\int\n^\eta|u|^2|\na u|^2dx &\le C\|\n^{\eta/2}u\|_{L^{8/\ti\eta}}^{2}\|\na u\|_{L^{8/(4-\ti\eta)}}^{2} \\&\le C(\eta)\psi^{4+2\eta}\left(1+\|\na^2 u\|_{L^2}^{\ti\eta/2} \right) \\&\le C(\ve,\eta) \psi^{\alpha(\eta)}+\ve \psi^{-2}\|\na^2 u\|_{L^2}^2.\ea\ee

Then, noticing that $\lambda=b\n^{\beta}$ satisfies
\be\ba\la{bb1} \lambda_{t} +\div( \lambda u)+ (\beta-1)  \lambda\div u=0,\ea\ee we obtain after using \eqref{cc5} and  \eqref{cc6} that
\be\la{bv1}\ba  \int \lm_t (\div u)^2dx & \le C\int \lm |u||\na u||\na^2u|dx+C\int \lm|\na u|^3dx\\& \le C(\ve)\psi\int \lm^2 |u|^2|\na u|^2 dx+\ve\psi^{-1}\|\na^2 u\|_{L^2}^2+C\psi^\beta\|\na u\|_{L^3}^3 \\& \le  C(\ve)\psi^\alpha  +C \ve\psi^{-1}\|\na^2 u\|_{L^2}^2 ,\ea\ee
where (and what follows)  we use  $\alpha=\alpha(\beta,\ga,q)>1$  to  denote   a  generic   constant   depending only  on $\beta,\ga,$ and $ q$, which may be different from line to line.

  Finally,
since $P $ satisfies
\be\la{bv2}  P_{t} +\div(P  u)+ ( \ga-1)P \div u=0,\ee
we deduce from \eqref{3.a2},  \eqref{cc6}, and the Sobolev inequality that
\be\la{bv3}\ba &2\int P\div u_tdx\\ &=2\frac{d}{dt}\int P\div udx-2\int Pu\cdot\na\div udx+2(\ga-1)\int P(\div u)^2dx   \\ &\le 2\frac{d}{dt}\int P\div udx+ \ve\psi^{-1}\|\na^2 u\|_{L^2}^2+C(\ve) \psi^\alpha.\ea\ee

Putting \eqref{cc5}, \eqref{bv1}, and \eqref{bv3} into \eqref{3r1} gives
\be\la{3.8-1}\ba& \frac{d}{dt}\int\left((2\mu+\lm)(\div u)^2+\o^2+R^{-1}|u|^2- 2P\div u\right)dx+\int \n|  u_t|^2dx\\ &\le C \ve\psi^{-1}\|\na^2 u\|_{L^2}^2+C(\ve) \psi^\alpha.\ea\ee

To estimate the first term on the right-hand side of \eqref{3.8-1}, it follows from \eqref{lp} and \eqref{2.9-1} that for $p\in [2,q],$
 \be\la{r3.22}   \ba  \|\na^2u\|_{L^p}\le &C\|\na \o\|_{L^p}+C\|\na \div u\|_{L^p}
\\ \le &C \left(\| \na  \o\|_{L^p} +  \|\na\left((2\mu+\lambda)\div u\right)\|_{L^p}+  \||\na \lambda|\div u\|_{L^{p}}  \right)\\
\le &C \left(\|\n \dot u\|_{L^p}  + \|\na P\|_{L^p}+R^{-1}\|u\|_{L^p} +  \||\na \lambda|\div u\|_{L^{p}}\right)  ,\ea\ee
which together with  \eqref{cc6}  and \eqref{cc5} leads to
  \be \la{cca}\ba   \|\na^2u\|_{L^2}\le &  C \psi^{1/2} \|\sqrt{\n} \dot u\|_{L^2}    +C\psi^{\alpha}+C\|\na \lambda\|_{L^{q}}  \|\na u\|_{L^{2 q/(q-2) }} \\ \le  & C \psi^{1/2} \|\sqrt{\n}   u_t\|_{L^2}    +C\psi^{\alpha}+\frac12 \|\na^2 u\|_{L^2}    .\ea\ee
Putting \eqref{cca} into \eqref{3.8-1},  integrating   the resulting inequality  over $(0,t),$ and choosing $\ve$ suitably small  yield that
 \bnn \ba   R^{-1}\|u\|_{L^2}^2+\|\na u\|^2_{L^{2}}
  +\int_{0}^{t}\int{\rho} |u_t|^2  dx ds
 &\le C+C \|P\|_{L^{2}}^{2}+C\int_{0}^{t}\psi^\alpha ds \\
 & \le C+C\int_{0}^{t} \psi^{\alpha} ds,
\ea
\enn where in the second inequality we have used
  \bnn\ba  \|P\|_{L^{2}}^{2}&\le \|P(\n_{0})\|_{L^{2}}^{2}+C\int_{0}^{t} \|P\|_{L^1}^{1/2} \|P\|_{L^\infty}^{3/2}\|\na u\|_{L^2}ds
  \le C+C\int_{0}^{t} \psi^{\alpha} ds,
 \ea\enn due to  \eqref{bv2}.
    The proof of  Lemma \ref{l3.0} is finished.

\begin{lemma}\la{l3.2}
 Let   $(\n,u)$ and $T_1$  be as   in Lemma \ref{l3.0}. Then for all $t\in (0, T_1],$
 \be\ba\la{li-1a}
 \sup_{0\le s\le t} s\int \n |u_{t}|^2dx+\int_{0}^{t}s\int \left( |\na u_t|^2 +R^{-1}|u_t|^2\right) dxds \le C   \exp\left\{C \int_0^t\psi^\alpha ds\right\}.\ea \ee
\end{lemma}

{\it Proof.}  Differentiating $\eqref{b2}_2$ with respect to $t$ gives
\be\la{zb1}\ba &\n u_{tt}+\n u\cdot \na u_t-\na ((2\mu+\lm)\div u_t)-\mu\na^\perp \o_t+R^{-1}u_t\\ &=-\n_t(u_t+u\cdot\na u)-\n u_t\cdot\na u+\na (\lm_t\div u)-\na P_t.\ea\ee
 Multiplying \eqref{zb1} by $u_t$ and integrating the resulting equation over $\O,$ we obtain after using  $\eqref{b2}_1$ that\be\ba  \la{na8}&\frac{1}{2}\frac{d}{dt} \int \n |u_t|^2dx+\int \left((2\mu+\lm)(\div u_t)^2+\mu \o_t^2+R^{-1}|u_t|^2\right)dx\\
 &=-2\int \n u \cdot \na  u_t\cdot u_tdx  -\int \n u \cdot\na (u\cdot\na u\cdot u_t)dx\\
  &\quad-\int \n u_t \cdot\na u \cdot  u_tdx
-\int \lm_t\div u \div u_{t}dx+\int P_{t}\div u_{t} dx\\
 &\le C\int  \n |u||u_{t}| \left(|\na  u_t|+|\na u|^{2}+|u||\na^{2}u|\right)dx +C\int \n |u|^{2}|\na u ||\na u_{t}|dx \\
 &\quad+C\int \n |u_t|^{2}|\na u |dx +C\int |\lm_t||\div u| |\div u_{t}|dx+C\int |P_{t}||\div u_{t}|dx.
  \ea\ee

  We estimate each term  on the right-hand side of  \eqref{na8} as follows:

First,
it follows from  \eqref{3.2}, \eqref{3.1-c},   \eqref {3.v2},
\eqref{3.a2}, and \eqref{cc6}  that for $\ve\in (0,1),$
 \be\la{na2}\ba  &\int  \n |u||u_{t}| \left(|\na  u_t|+|\na u|^{2}+|u||\na^{2}u|\right) dx\\
 & \le C \|\n^{1/2} u\|_{L^{6}}\|\n^{1/2} u_{t}\|_{L^{2}}^{1/2} \|\n^{1/2} u_{t}\|_{L^{6}}^{1/2}\left(\| \na u_{t}\|_{L^{2}}+\| \na u\|_{L^{4}}^{2} \right) \\
 &\quad +C\|\n^{1/4}  u \|_{L^{12}}^{2}\|\n^{1/2} u_{t}\|_{L^{2}}^{1/2} \|\n^{1/2} u_{t}\|_{L^{6}}^{1/2} \| \na^{2} u \|_{L^{2}}  \\
 & \le C\psi^{\alpha} \|\n^{1/2} u_{t}\|_{L^{2}}^{1/2}\left(\|\n^{1/2} u_{t}\|_{L^{2}} +\| \na u_{t}\|_{L^{2}}\right)^{1/2}\left(\| \na u_{t}\|_{L^{2}}+  \| \na^{2} u \|_{L^{2}}+\psi\right)\\
 &\le  \ve\| \na u_{t}\|_{L^{2}}^{2}+C(\ve)\psi^{\alpha} \left(\| \na^{2} u \|_{L^{2}}^{2} +  \|\n^{1/2} u_{t}\|_{L^{2}}^{2}+1\right).
 \ea\ee

 Next, Holder's inequality together with \eqref{3.a2}   and \eqref{cc6} yields that
 \be \ba \la{5.ap3}  \int \n |u|^{2}|\na u ||\na u_{t} |dx&\le C \|\n^{1/2} u\|_{L^{8}}^{2}\|\na u\|_{L^{4}} \| \na u_{t}\|_{L^{2}}\\
&\le  \ve \| \na u_{t}\|_{L^{2}}^{2}+  C(\ve)\left(\psi^{\alpha}+\| \na^{2} u \|_{L^{2}}^{2}\right).\ea\ee

Then, Holder's inequality and \eqref{3.v2}  lead to
 \be \ba \la{5.a3}  \int \n |u_t|^{2}|\na u |dx&\le   \| \na u\|_{L^{2}}
   \|\n^{1/2} u_{t}\|_{L^{6}}^{3/2}\|\n^{1/2} u_{t}\|_{L^{2}}^{1/2} \\
  &\le  \ve \| \na u_{t}\|_{L^{2}}^{2} + C(\ve)\psi^{\alpha}  \|\n^{1/2} u_{t}\|_{L^{2}}^{2}.\ea\ee

Next, we use \eqref{bb1} and  \eqref{3.a2}  to  get \be\la{na3}\ba &\int |\lm_t||\div u| |\div u_{t}|dx\\
&\le C\int \left(\lm(\div u)^{2} |\div u_{t}|+|\na \lm| |u| |\div u||\div u_{t}|\right)dx\\
&   \le \frac{1}{2}\int  \lm(\div u_{t})^{2}dx +C\psi^{\beta} \|\na u\|_{L^{4}}^{4}\\
&\quad +C\|\bar{x}^{\te_{0}}\na \lm \|_{L^{q}}\|  u \bar{x}^{- \te_{0}}\|_{L^{4q/((q-2) \te_{0})}}\|\na u \|_{L^{4q/((q-2)(2- \te_{0}) )}}\|\na u_{t}\|_{L^{2}} \\
&\le  \frac{1}{2}\int  \lm(\div u_{t})^{2}dx+ \ve \|\na u_{t}\|_{L^{2}}^{2} + C(\ve ) \psi^{\alpha} \| \na^{2} u \|_{L^{2}}^{2}+ C(\ve ) \psi^{\alpha} .\ea\ee

Finally, it follows from  \eqref{bv2} and  \eqref{3.a2}   that
 \be \la{na1}\ba & \int |P_{t}||\div u_{t}|dx\\
 &\le C  \int \left(P|\div u|+ |\na P||u|\right)|\div u_{t}|dx\\
  &\le C\left(\|P\|_{L^{\infty}}\|\na u\|_{L^{2}} + \|\n^{\ga-1}\|_{L^{\infty}}\|\bar{x}^{a}\na\n \|_{L^{q}}\|\bar{x}^{-a}u\|_{L^{2q/(q-2)}} \right)\|\na u_{t}\|_{L^{2}} \\  &\le \ve \|\na u_{t}\|_{L^{2}}^{2} + C(\ve) \psi^{\alpha},\ea\ee
 where in the last inequality we have used  \eqref{3.a2}.

 Substituting \eqref{na2}-\eqref{na1} into \eqref{na8} and choosing $\ve$ suitably  small lead to
\be\ba\la{a4.6} & \frac{d}{dt} \int \n |u_t|^2dx+\int \left((2\mu+\lm)(\div u_t)^2+ \mu\o_t^2+R^{-1}|u_t|^2\right)dx\\
  &  \le   C\psi^{\alpha}\left(1+\| \n^{1/2}u_{t}\|_{L^{2}}^{2}+\| \na^{2}u\|_{L^{2}}^{2}\right)\\
 &  \le   C \psi^{\alpha}  \|\n^{1/2} u_{t}\|_{L^{2}}^{2}+C\psi^{\alpha},
  \ea\ee
 where  in the last inequality we have used \eqref{cca}.
Multiplying  \eqref{a4.6}  by $t , $  we obtain \eqref{li-1a} after
using Gronwall's  inequality and \eqref{3.1}. The proof of  Lemma \ref{l3.2} is completed.

\begin{lemma}\la{l3.4}
 Let   $(\n,u)$ and $T_1$  be as   in Lemma \ref{l3.0}. Then for all $t\in (0, T_1],$
  \be\la{b1a} \ba &\sup\limits_{0\le s\le t}\left( \|\n \bar x^a\|_{L^1\cap H^1 \cap  W^{1,q}}+\|\lm(\n)\|_{L^2}+\|\bar{x}^{\te_{0}}\na \lm(\n)\|_{L^2\cap L^q} \right)\\&  \le \exp\left\{C\exp\left\{C\int_{0}^{t} \psi^{\alpha} ds\right\}\right\} . \ea \ee
\end{lemma}

 \emph{Proof}. First, \eqref{bb1} gives
  \be \la{8.1}(\|\lm(\n)\|_{L^2})_t\le C\|\na u\|_{L^\infty}\| \lm(\n)\|_{L^2}\le C\|\na u\|_{W^{1,q}}\| \lm(\n)\|_{L^2}.\ee

Next, multiplying $\eqref{b2}_1$ by $\bar x^\al $ and integrating the resulting equality over $\O,$ we obtain after integration by parts and using \eqref{3.1-c} that
\bnn \ba  \frac{d}{dt}\int\n \bar x^\al  dx &\le C\int \n |u|\bar x^{\al -1}\log^{1+\eta_0}(e+|x|^2)dx\\ &\le C\left(\int\n \bar x^{2\al -2} \log^{2(1+\eta_0)}(e+|x|^2)dx\right)^{1/2}\left(\int\n u^2 dx\right)^{1/2}\\ &\le C\left(\int\n \bar x^{ \al  } dx\right)^{1/2},\ea\enn
which   gives
\be\la{o3.7} \sup_{0\le t\le T_1}\int\n \bar x^\al  dx\le C .\ee

 Next, it follows from the Sobolev inequality and \eqref{3.a2} that
for $0<\de<1,$   \be\ba\la{3.22}
  \norm[L^{\infty}]{u\bar x^{-\de}}
&\le C(\de)\left(\norm[L^{4/\de}]{u\bar x^{-\de}}+ \norm[L^{3}]{\na (u\bar x^{-\de})}\right) \\
& \le C(\de)\left(\norm[L^{4/\de}]{u\bar x^{-\de}}+ \norm[L^{3}]{\na u}+\norm[L^{4/\de}]{u\bar x^{-\de}} \| \bar {x}^{-1}\na\bar{x}  \|_{L^{12/(4-3\de)}} \right) \\
& \le C(\de) \left(\psi^\alpha+\norm[L^2]{\na^2 u} \right).\ea\ee
 One derives from $\eqref{n1}_1$ that $ w\triangleq\n\bar x^a$ satisfies\bnn  w_t+u\cdot\na w-a wu\cdot\na \log \bar x+w\div u=0,\enn
which together with \eqref{3.22} gives that for   $p\in [2,q]$
\be\la{7.z1}\ba (\|\na w\|_{L^p} )_t& \le C(1+\|\na u\|_{L^\infty}+\|u\cdot \na \log \bar x\|_{L^\infty}) \|\na w\|_{L^p} \\&\quad +C\left( \||\na u||\na\log \bar x|\|_{L^p}+\||  u||\na^2\log \bar x|\|_{L^p}+\| \na^2 u \|_{L^p}\right)\|w\|_{L^\infty} \\& \le C(\psi^\alpha  +\|\na u\|_{L^2\cap W^{1,q}})  \|\na w\|_{L^p} \\&\quad+C\left(\|\na u\|_{L^p}+\|u\bar x^{-1/4}\|_{L^\infty}\|\bar x^{-3/2}\|_{L^p}+\|\na^2 u\|_{L^p}\right) \|w\|_{L^\infty}\\& \le C(\psi^\alpha  +\|\na^2 u\|_{L^2\cap L^{q}})(1+ \|\na w\|_{L^p}+\|\na w\|_{L^q}),\ea\ee where in the last inequality we have used \eqref{o3.7}.
Similarly, one obtains from \eqref{bb1} that
\be\la{7.z2}\ba  (\|\na (\bar x^{\theta_0}\lm)\|_{L^2 \cap L^q}  )_t  \le C(\psi^\alpha  +\|\na^2 u\|_{L^2\cap L^{q}})(1+\|\na (\bar x^{\theta_0}\lm)\|_{L^2 \cap L^q}   ).  \ea\ee

 Next, we claim that  \be\ba  \la{bb6-1}  \int_{0}^{t}\left(\|\na^2u\|_{L^2\cap L^q}^{(q+1)/q} +t\|\na^2u\|_{L^2\cap L^q}^2\right)dt
  \le  C   \exp\left\{C \int_0^t\psi^\alpha ds\right\},\ea\ee
 which together with \eqref{7.z1}, \eqref{7.z2},
  \eqref{o3.7}, \eqref{8.1},      and the Gronwall  inequality yields that
 \be\la{b2a} \ba &\sup\limits_{0\le s\le t}\left(\| \bar x^a \n\|_{  L^1\cap H^1\cap W^{1,q}}+\|\lm(\n)\|_{L^2}+\|\na(\bar{x}^{\te_{0}} \lm)\|_{L^2\cap L^q} \right)\\&\le \exp\left\{C\exp\left\{C\int_{0}^{t} \psi^{\alpha} ds\right\}\right\} . \ea \ee One thus directly obtains \eqref{b1a} from this and the following simple fact:
 \bnn \|\bar{x}^{\te_{0}}\na \lm\|_{L^2\cap L^q} \le \|\na(\bar{x}^{\te_{0}} \lm)\|_{L^2\cap L^q}+C\| \bar x^a \n\|^\beta_{  W^{1,q}},\enn  due to  $\te_0<a\beta.$

Finally, to finish  the proof of    Lemma \ref{l3.4}, it only remains to prove \eqref{bb6-1}. In fact, on the one hand, it follows  from \eqref{cca}, \eqref{3.1}, and \eqref{li-1a}  that  \be\ba   \la{bb15}&\int_0^{t} \left(\|\na^2u\|_{L^2}^{5/3}+s\|\na^2u\|_{L^2}^2 \right)ds\\&
 \le  C \int_0^{t}\left( \|\sqrt{\n }  u_t\|_{L^2}^{2}  + \psi^{\alpha}  \right)ds+C   \exp\left\{C \int_0^t\psi^\alpha ds\right\}\int_0^t\psi^\alpha ds
  \\& \le C   \exp\left\{C \int_0^t\psi^\alpha ds\right\}.\ea\ee
On the other hand, choosing $p=q$ in \eqref{r3.22} gives \be \la{3.1a5} \ba  \|\na^2u\|_{L^q}\le  &C \left(\|\n \dot u\|_{L^q}  + \|\na P\|_{L^q}+R^{-1}\|u\|_{L^q} + \|\na \lambda\|_{L^{q}}  \|\na u\|_{L^\infty}\right) \\
\le &C \left(\|\n \dot u\|_{L^q}   + \psi^\alpha +  \psi^\alpha\|\na^2 u\|_{L^q}^{q/ (2q-2)}\right)\\
\le & \frac{1}{2}\| \na^{2} u \|_{L^{q}}  +C\psi^{\alpha} +C   \|\n \dot u\|_{L^q}.\ea\ee
By \eqref{3.v2}, \eqref{3.a2}, and \eqref{cc6},  the last term on the right-hand side of \eqref{3.1a5} can be estimated  as follows: \bnn\ba \|\n \dot u\|_{L^q}&\le \|\n   u_t\|_{L^q}+\|\n u\cdot \na u\|_{L^q}\\ &\le \|\n   u_t\|_{L^2}^{2(q-1)/(q^2-2)}\|\n   u_t\|_{L^{q^2}}^{(q^2-2q)/(q^2-2)}+\|\n u\|_{L^{2q}} \|\na u\|_{L^{2q}}\\ &\le C\psi^\alpha \left( \|\n^{1/2}   u_t\|_{L^2}^{2(q-1)/(q^2-2)}\|\na   u_t\|_{L^{2 }}^{(q^2-2q)/(q^2-2)} +\|\n^{1/2}   u_t\|_{L^2}\right) \\&\quad +C\psi^\alpha \left(1+ \|\na^2 u\|^{1-1/q}_{L^2 }\right).\ea\enn
This combined with \eqref{bb15},   \eqref{li-1a},  and \eqref{3.1} yields that
 \be \la{3.1a6}  \ba  \int_0^t\|\n \dot u\|_{L^q}^{\frac{q+1}{q}}dt &\le  C  \int_0^t\psi^\alpha t^{-\frac{q+1}{2q}} \left(t\|\n^{1/2}u_t\|_{L^2}^2\right)^{\frac{q^2-1}{ q(q^2-2)}} \left(t\|\na u_t\|^2_{L^2}\right)^{\frac{(q-2)(q+1)}{2 (q^2-2)}} dt\\&\quad+C\int_0^t\|\n^{1/2}u_t\|_{L^2}^2dt +C   \exp\left\{C \int_0^t\psi^\alpha ds\right\} \\&\le  C   \exp\left\{C \int_0^t\psi^\alpha ds\right\} \int_0^t\left(\psi^\alpha+t^{-\frac{q^3+q^2-2q-1}{q^3+q^2-2q}}+t\|\na u_t\|_{L^2}^2\right) dt \\&\quad+C   \exp\left\{C \int_0^t\psi^\alpha ds\right\}  \\&\le  C   \exp\left\{C \int_0^t\psi^\alpha ds\right\}, \ea\ee
 and that \be\la{3.1a7}  \int_0^t t\|\n \dot u\|_{L^q}^2dt \le  C   \exp\left\{C \int_0^t\psi^\alpha ds\right\}.\ee One thus obtains \eqref{bb6-1} from  \eqref{bb15}-\eqref{3.1a7}   and finishes the proof of    Lemma \ref{l3.4}.

   Now,
  Proposition \ref{pro} is a direct consequence of  Lemmas \ref{l3.0}-\ref{l3.4}.

\emph{Proof of Proposition \ref{pro}}. It follows from \eqref{3.1-c}, \eqref{3.1}, and \eqref{b1a} that
\bnn\ba \psi(t)
&\le   \exp\left\{C\exp\left\{C\int_{0}^{t} \psi^{\alpha} ds\right\}\right\}.\ea\enn Standard arguments thus yield  that for $M\triangleq e^{Ce}$ and $T_0\triangleq \min\{T_1,(CM^\alpha)^{-1}\} ,$
\bnn \sup\limits_{0\le t\le T_0}\psi(t)\le M,\enn which together with  \eqref{bb6-1},    \eqref{cca}, and \eqref{3.1}  gives  \eqref{o1}.
The proof of Proposition \ref{pro} is thus completed.

\section{A priori estimates (II)}

In this section, in addition to $\mu,  \beta,  \gamma, b, q, a, \eta_0, \te_0,  N_0, $ and $ E_0, $  the generic positive
constant $C $ may depend   on   $\de_0, \|\na^2u_0\|_{L^2},$ $\|\na^2  \n_0\|_{L^q}, $ $\|\na^2  \lm(\n_0)\|_{L^q}, $ $\|\na^2  P(\n_0)\|_{L^q}, $  $\| \bar x^{\de_0}\na^2  \rho_0 \|_{L^2},$ $\| \bar x^{\de_0}\na^2  \lm(\rho_0) \|_{L^2},$   $\| \bar x^{\de_0}\na^2  P(\rho_0 )\|_{L^2},$  and $\|\ti g\|_{L^2},$
 where  \be\la{gq1} \ti g\triangleq\n_0^{-1/2}(-\mu\Delta u_0-\na((\mu+\lm(\n_0))\div u_0 )+\na P(\n_0)+R^{-1}u_0).\ee

\begin{lemma}\label{lem4.5}  It holds that
\begin{equation}\la{5.13a}\ba
 \sup_{0\leq t\leq T_0}\left(\|\bar x^{ \de_0}\na^2\n\|_{L^2}+\|\bar x^{ \de_0}\na^2\lm \|_{L^2}+\|\bar x^{ \de_0}\na^2 P \|_{L^2}\right)\le C.\ea
\end{equation}
\end{lemma}

\pf First, by virue of \eqref{2.1} and $\eqref{b2}_2,$  defining $$\n^{1/2}u_t(x,t=0)\triangleq- \ti g-\n_0^{1/2} u_0\cdot \na u_0,$$     integrating \eqref{a4.6} over $(0,T_0),$    and using \eqref{o1} and \eqref{3.1},   we obtain that
 \be\la{oo1} \sup\limits_{0\le t\le T_0}\int \n |u_t|^2dx+\int_0^{T_0} \left(\|\na u_t\|_{L^2}^2+R^{-1}\|u_t\|_{L^2}^2\right)dt\le C.\ee
  This combined with \eqref{cca} and \eqref{o1} gives \be\la{oo2} \sup\limits_{0\le t\le T_0}\|\na u\|_{H^1}\le C,\ee which together with \eqref{3.22} and \eqref{o1}   shows that for   $\de\in(0,1),$ \be\la{oo3} \|\n^\de u\|_{L^\infty}+\|\bar x^{-\de}u\|_{L^\infty}\le C(\de).\ee Direct calculations  yield  that for  $2\le r\le q$
\be \la{oo7} \|(\bar x^{(1+a)/2}+|u|) \rho_{t}\|_{L^{r}} +\|(1+|u|)\lm_{t}\|_{L^{r}}+\|(1+|u|)P_{t} \|_{L^{r}}\le C,\ee due to \eqref{o1},  $\eqref{b2}_1$, \eqref{bb1},   \eqref{bv2}, \eqref{oo2}, and \eqref{oo3}.
It follows from \eqref{3.v2}   and \eqref{oo1}-\eqref{oo3}   that
for   $\de\in (0,1] $ and   $s> 2/\de,$
\be\la{5.22}\ba \|\bar x^{-\de}u_t\|_{L^s}+\|\bar x^{-\de}\dot u\|_{L^s}&\le C\|\bar x^{-\de}  u_t\|_{L^s}+C\|\bar x^{-\de}  u \|_{L^\infty}\|\na  u \|_{L^s}\\& \le C(\de,s)+C(\de,s)\|\na u_t\|_{L^2}.\ea\ee

Next,
denoting  $v\triangleq  \bar x^{\de_0} g(\n)$ with $g(\n)=\n^p$ for   $ p\in[ \min\{\beta,1\},\max\{\beta,\ga\}],$  we get from \eqref{o1} that \be\la{o4}\|v\|_{L^\infty}+\|\na v\|_{L^2\cap L^q} \le C,\ee due to $\de_0\le \te_0\le \min\{1, \beta \} .$ It follows from   $(\ref{b2})_1$ that \be\la{a4.r3} g(\n)_t+u\cdot \na g(\n)+pg(\n)\div u=0,\ee which gives \bnn v_t+u\cdot \na v-\de_0 v u\cdot \na \log \bar x+pv\div u=0 .\enn Thus,
direct calculations yield   that
\be \la{oj1}\ba  (\|\na^2 v\|_{L^2})_t &\le C(1+\|\na u\|_{L^\infty}+\|u\cdot \na\log \bar x\|_{L^\infty})\|\na^2v\|_{L^2} +C\||\na^2u||\na v|\|_{L^2}\\&\quad +C\||\na v||\na u||\na \log\bar x|\|_{L^2}+C\||\na v||u| |\na^2 \log\bar x|\|_{L^2}\\&\quad +C \|v\|_{L^\infty}  \left(\|\na^2 (u\cdot\na \log \bar x)\|_{L^2}+\|\na^3u\|_{L^2}\right)\\&\le C(1+\|\na u\|_{L^\infty} )\|\na^2v\|_{L^2} +C\| \na^2u\|_{L^{2q/(q-2)}}\|\na v \|_{L^q}\\&\quad+C\| \na v\|_{L^2}\|\na u\|_{L^\infty} +C\| \na v\|_{L^2}\||u| |\na^2 \log\bar x |\|_{L^\infty} \\&\quad+C   \|\na^2  u \|_{L^2}+C\|\na   u \|_{L^2} +C\|   u|\na^3 \log \bar x|\|_{L^2}+C\|\na^3u\|_{L^2}\\&\le C(1+\|\na  u\|_{L^\infty} ) \|\na^2v\|_{L^2}  +C+C \|\na^3u\|_{L^2} ,\ea\ee
where in the second and third inequalities we have used  \eqref{oo3} and \eqref{o4}. We use \eqref{lp}, \eqref{oo2},  and \eqref{2.9-1}  to estimate the last term on the right-hand side of \eqref{oj1}  as follows:
\be \label{4.d25}\ba
\|\nabla^3 u\|_{L^2}   &\le  C \|\na(\na^\perp \o)\|_{L^2}+C\|\na(\na \div u)\|_{L^2} +C\\ &\le  C\left(\|\na^2 \o \|_{L^2}+\|\na^2 P \|_{L^2}+\|\na^2 F\|_{L^2}+\||\na^2u||\na \lm |\|_{L^2}\right)\\&\quad+C\||\na u||\na^2 \lm| \|_{L^2}+C\\ &\le  C\left(\|\na (\n\dot u)\|_{L^2}+\|\bar x^{\de_0}\na^2 P \|_{L^2} +  \|\na \lm \|_{L^q}\|\na^2u\|_{L^{2q/(q-2)}}\right)\\&\quad+C \|\bar x^{\de_0}\na^2 \lm \|_{L^2}\|\na u\|_{L^\infty}+C\\ &\le  C \| \na u_t \|_{L^2}+C\|\bar x^{\de_0}\na^2 P \|_{L^2} + \frac12 \|\na^3u\|_{L^2}\\&\quad+C \|\bar x^{\de_0}\na^2 \lm \|_{L^2}\|\na u\|_{L^\infty}+C ,
\ea\ee where in the last inequality we have used the following simple fact:
\bnn\ba   \|\na(\n \dot u)\|_{L^2}&\le C\|\n \bar x^a\|_{L^6} \|\bar x^{-a}\dot u \|_{L^3} +C\|\bar x^a\na\n\|_{L^q}\|\bar x^{-a}\dot u\|_{L^{2q/(q-2)}}\\&\le C\|\na u_t\|_{L^2}+C,\ea\enn due to \eqref{o1},  \eqref{5.22}, and \eqref{oo2}.
  Noticing that  \eqref{o1} leads to
\be\la{qq}\ba &
 \|\bar x^{ \de_0}\na^2\n\|_{L^2}+\|\bar x^{ \de_0}\na^2\lm \|_{L^2}+\|\bar x^{ \de_0}\na^2 P \|_{L^2} \\&\le C\|\na^2 (\bar x^{\de_0}\n )\|_{L^2}+ C\|\na^2 (\bar x^{\de_0}\lm )\|_{L^2}+C\|\na^2 (\bar x^{\de_0}P)\|_{L^2}+C,\ea\ee we substitute   \eqref{4.d25} into \eqref{oj1} to get
\bnn\ba&
\frac{d}{dt}\left(\|\na^2 (\bar x^{\de_0}\n )\|_{L^2}+ \|\na^2 (\bar x^{\de_0}\lm )\|_{L^2}+\|\na^2 (\bar x^{\de_0}P)\|_{L^2}\right)\\ &\leq  C (1+\|\na^2u\|_{L^q})\left(\|\na^2 (\bar x^{\de_0}\n )\|_{L^2}+ \|\na^2 (\bar x^{\de_0}\lm )\|_{L^2}+\|\na^2 (\bar x^{\de_0}P)\|_{L^2}\right)\\& \quad+C\|\na u_t\|_{L^2}+C,
\ea\enn
which, along   with Gronwall's inequality,  \eqref{o1},  and   \eqref{oo1}, yields that \be\la{5.13b} \sup\limits_{0\le t\le T_0}\left(\|\na^2 (\bar x^{\de_0}\n )\|_{L^2}+ \|\na^2 (\bar x^{\de_0}\lm )\|_{L^2}+\|\na^2 (\bar x^{\de_0}P)\|_{L^2}\right)\le C. \ee This  combined with  \eqref{qq}    gives \eqref{5.13a} and finishes the proof of Lemma \ref{lem4.5}.

 \begin{lemma} \la{oiq1} It holds that
\be\ba \la{oo8}\sup\limits_{0\le t\le T_0}t\left(\|\na u_t\|_{L^{2}}^2 +R^{-1}\|u_t\|_{L^{2}}^2\right) +\int_{0}^{T_0}t\left(\|\n^{1/2}  u_{tt}\|_{L^{2}}^2+\|\na^2  u_{t}\|_{L^{2}}^2\right)dt \le C .\ea\ee
 \end{lemma}

{\it Proof.} Multiplying \eqref{zb1} by $u_{tt}$ and integrating the resulting equation over $\O$ lead to
\be\ba \la{oo4}&\frac{1}{2}\frac{d}{dt}\int\left((2\mu+\lm)(\div u_t)^2+\mu \o_t^2+R^{-1}|u_t|^2\right)dx+\int\n |u_{tt}|^2dx\\
&=-\int\left(2 \n u\cdot\na u_t\cdot u_{tt}+ \n u_t \cdot\na u\cdot  u_{tt}\right)dx -\int \n u\cdot \na(u \cdot\na u)  \cdot u_{tt}dx\\
&\quad-\int \n u\cdot  \na   u_{tt}\cdot \dot u  dx +\frac12\int\lm_t(\div u_t)^2dx
-\int\lm_t\div u\div u_{tt}dx\\&\quad+\int P_t\div u_{tt}dx.\ea\ee

We   estimate each term on the right-hand side of \eqref{oo4} as follows:

  First, it follows from   \eqref{o1}, \eqref{oo1}-\eqref{oo3}, and \eqref{5.22} that\be\la{6.15}\ba & \left|\int\left(2 \n u\cdot\na u_t\cdot u_{tt}+ \n u_t \cdot\na u\cdot  u_{tt}\right)dx \right|+\left|\int \n u\cdot \na(u \cdot\na u)  \cdot u_{tt}dx\right|\\
  &\le \ep \int\n |u_{tt}|^2dx +C(\ep)\left(\|\n^{1/2} u\|_{L^{\infty}}^{2}\|\na u_{t}\|_{L^{2}}^{2}+\|\n^{1/2} u_{t}\|_{L^{4}}^{2}\|\na u\|_{L^{4}}^{2}\right)\\&\quad +C(\ep)\int\left( \n |u|^{2}|\na u|^{4}+ \n |u|^4 |\na^{2} u|^{2}\right)dx \\
    &\le \ep \int\n |u_{tt}|^2dx +C(\ep)\left( \|\na u_{t}\|_{L^{2}}^{2}+1\right).
   \ea \ee

Next, direct calculations give
\be\la{oo16}\ba -\int \n u\cdot\na u_{tt}\cdot \dot udx
&=-\frac{d}{dt}\int \n u\cdot\na u_{t}\cdot\dot udx +\int  (\n u )_t\cdot \na  u_{t}\cdot \dot u    dx\\&\quad+ \int \n u\cdot\na u_t\cdot\left(u_{tt}+ u_{t} \cdot\na u+ u \cdot\na u_{t} \right)  dx .
\ea\ee
On the one hand,  it follows from    \eqref{o1}  and \eqref{oo3}-\eqref{5.22} that
\be\la{oo14}\ba  \left| \int  (\n u)_t\cdot \na  u_{t}\cdot \dot u  dx\right|   & \le C\|\n\bar x^a\|_{L^\infty}\|\bar x^{-a/2} u_t \|_{L^4}\|\bar x^{-a/2} \dot u \|_{L^4}\|\na u_t\|_{L^2}\\&  \quad+C \|\bar x^{(a+1)/2}\n_t\|_{L^2}  \|u\bar x^{-1/2}\|_{L^\infty}\| \bar x^{-a/2}\dot u\|_{L^4} \|\na u_t\|_{L^4}\\& \le  C(\de)\|\na u_t\|^4_{L^2}+C(\de)+\de\|\na^2 u_t\|_{L^2}^2 .\ea\ee
On the other hand, Cauchy's inequality   and   \eqref{oo1}-\eqref{5.22} lead to
\be\la{oo15}\ba  &\left|    \int \n u\cdot\na u_t\cdot \left(u_{tt}+ u_{t} \cdot\na u+ u \cdot\na u_{t} \right) dx\right|\\&\le \ve\int \n |u_{tt}|^2dx+C(\ve)\|\n^{1/2}  u\|_{L^\infty}^2\int|\na u_t|^2dx\\&\quad+C\|\bar x^{-a/2}u\|_{L^\infty} \|\na u_t\|_{L^2} \|\bar x^{-a/2}u_t\|_{L^4}\|\na u\|_{L^4} +C \|\bar x^{-1/2}u\|_{L^\infty}^2\|\na u_t\|_{L^2}^2\\&\le \ve\int \n |u_{tt}|^2dx+C(\ve)\int |\na u_t|^2dx+C(\ve).\ea\ee
Putting \eqref{oo14} and \eqref{oo15} into \eqref{oo16} thus shows
\be\la{5.za61}\ba  -\int \n u\cdot\na u_{tt}\cdot \dot u dx
&\le-\frac{d}{dt}\int \n u\cdot\na u_{t}\cdot\dot u dx+\ve\int \n |u_{tt}|^2dx\\&\quad+\de\|\na^2 u_t\|_{L^2}^2+C(\ve,\de)\|\na u_t\|_{L^2}^4+C(\ve,\de).
\ea\ee

Next, the Sobolev inequality and \eqref{oo7} ensure
 \be\ba\la{6.15-a}  \int\lm_t(\div u_t)^2dx
 \le C \|\lm_{t}\|_{L^{2}} \|\na  u_t\|_{L^{4}}^2\le C(\de)  \|\na  u_t\|_{L^{2}}^2+ \de \|\na^{2}  u_t\|_{L^{2}}^2.
\ea\ee

Then,    \eqref{bb1}  leads to
\be\ba\la{3.2t}
&-\int\lm_t\div u\div u_{tt}dx\\
&=-\frac{d}{dt}\int\lm_t\div u\div u_{t}dx+\int\lm_t(\div u_t)^2dx\\ &\quad -(\beta-1)\int (\lm\div u)_t\div u\div u_tdx+\int(\lm u)_{t}\cdot\na(\div u\div u_{t})dx.
\ea\ee
It follows from \eqref{oo7}, \eqref{5.22}, and \eqref{oo2} that
\be\la{ko1}\ba& \left| \int (\lm\div u)_t\div u\div u_tdx \right|\\& \le C \left(\|\lm_t\|_{L^2}\|\na u\|_{L^8}^2\|\na u_t\|_{L^4}+\|\lm\|_{L^\infty}\|\na u \|_{L^2}\|\na u_t\|_{L^4}^2\right)\\&\le \de \|\na^2u_t\|_{L^2}^2+C(\de) \|\na u_t\|_{L^2}^2+C(\de)  ,\ea\ee and that
\be\la{ko2}\ba &\left| \int(\lm u)_{t}\cdot\na(\div u\div u_{t})dx\right|\\&\le C\|\lm_tu\|_{L^q} \left(\|\na u\|_{L^{2q/(q-2)}}\|\na^2 u_t\|_{L^2}+\|\na u_t\|_{L^{2q/(q-2)}}\|\na^2 u\|_{L^2}\right)\\&\quad+ C\|\lm u_t\|_{L^{4/(a\ti\beta)}} \left(\|\na u\|_{L^{4/(2-a\ti\beta)}}\|\na^2 u_t\|_{L^2}+\|\na u_t\|_{L^{4/(2-a\ti\beta)}}\|\na^2 u\|_{L^2}\right)\\&\le \de \|\na^2u_t\|_{L^2}^2+C(\de)(1+\|\na u_t\|_{L^2}^2),\ea\ee
with $\ti\beta=\min\{1,\beta\}.$  Putting \eqref{6.15-a}, \eqref{ko1}, and \eqref{ko2} into \eqref{3.2t} gives
 \be\ba\la{3.2s}
&-\int\lm_t\div u\div u_{tt}dx
\\&\le-\frac{d}{dt}\int\lm_t\div u\div u_{t}dx+C\de \|\na^2u_t\|_{L^2}^2 +C(\de)(1+\|\na u_t\|_{L^2}^2).
\ea\ee

 Finally, it follows from \eqref{bv2},  \eqref{oo7},  and \eqref{oo2} that
\be\ba\la{ss4-1}  &\int P_t\div u_{tt}dx\\
&=\frac{d}{dt}\int P_t\div u_{t}dx- \int (Pu)_{t}\cdot\na \div u_{t}dx +(\ga-1) \int (P\div u)_{t}  \div u_{t}dx\\
&\le \frac{d}{dt}\int P_t\div u_{t}dx+ C \left(\|P_tu\|_{L^2}+\|Pu_t\|_{L^2}\right)\|\na^2  u_{t}\|_{L^{2}} \\&\quad+C\|P_t\|_{L^2}\|\na u\|_{L^4}\|\na u_t\|_{L^4} +C\|P\|_{L^\infty} \|\na u_t\|_{L^2}^2\\
&\le \frac{d}{dt}\int P_t\div u_{t}dx+\de \|\na^2u_t\|_{L^2}^2+C(\de)(1+\|\na u_t\|_{L^2}^2).\ea \ee

Substituting \eqref{6.15},   \eqref{5.za61},  \eqref{6.15-a},  \eqref{3.2s}, and  \eqref{ss4-1}  into  \eqref{oo4} and choosing $\ve$ suitably small lead  to
 \be\la{oo11}\ba\Psi'(t)+\int\n|u_{tt}|^2dx\le C\de \|\na^2u_t\|_{L^2}^2+C(\de)\|\na u_t\|_{L^2}^4+C(\de),\ea\ee
where
\bnn\ba  \Psi(t)\triangleq &\int\left((2\mu+\lm)(\div u_t)^2+\mu \o_t^2+R^{-1}|u_t|^2\right)dx\\&
-2\int\left( P_t\div u_{t }- \lm_t\div u\div u_{t}-\n u\cdot\na u_{t}\cdot\dot u\right)dx  \ea\enn  satisfies
\be\la{o3}\ba   C_0(\mu)\|\na u_t\|_{L^2}^2+R^{-1}\|u_t\|_{L^2}^2-C\le \Psi(t)\le C\|\na u_t\|_{L^2}^2+R^{-1}\|u_t\|_{L^2}^2+C,\ea\ee
  due to the following simple fact:
  \bnn \ba &\left| \int\left( P_t\div u_{t }- \lm_t\div u\div u_{t}-\n u\cdot\na u_{t}\cdot\dot u\right)dx\right|\\ & \le C\|P_t\|_{L^2}\|\na u_t\|_{L^2}+C\|\lm_t\|_{L^q}\|\na u\|_{L^{2q/(q-2)}}\|\na u_t\|_{L^2}\\&\quad +C\|\n^{1/2}  u\|_{L^\infty}\|\na u_t\|_{L^2}\left(\|\n^{1/2}  u_t\|_{L^2}+\|\n^{1/2}u\cdot \na u\|_{L^2}\right)\\ &\le  \ve \|\na u_t\|_{L^2}^2+C(\ve),\ea\enn which comes from \eqref{oo1}-\eqref{oo7} and \eqref{cc5}.

Then, it remains to estimate the first term on the right-hand side of  \eqref{oo11}.   In fact, we obtain  from  \eqref{bb1},    \eqref{oo3}, and \eqref{5.13a}  that
 \be \la{o6}\ba\|\na\lm_t\|_{L^2}&\le C\|\na u\|_{L^\infty}\|\na\lm\|_{L^2}+C\|u\bar x^{-\de_0}\|_{L^\infty}\|\bar x^{\de_0}\na^2\lm\|_{L^2}+C\|\na^3u\|_{L^2}\\&\le C+C\|\na^3u\|_{L^2}\\&\le C+C\|\na u_t\|_{L^2},\ea\ee where in the last inequality we have used \be \la{5.b3} \|\na u\|_{H^2} \le C+C\|\na u_t\|_{L^2} ,\ee
  due to \eqref{4.d25},   \eqref{5.13b}, and \eqref{oo2}.
Similar to \eqref{o6}, we have
 \be\la{o7} \ba\|\na P_t\|_{L^2} \le C+C\|\na u_t\|_{L^2}.\ea\ee
Using the boundary condition $\eqref{b2}_3,$ we obtain   from \eqref{zb1} that
\bnn\ba & \| \na ((2\mu+\lm)\div u_t)\|_{L^2}^2+\mu^2\|\na^\perp \o_t\|_{L^2}^2 \\&\le \int \left| \na ((2\mu+\lm)\div u_t)+\mu\na^\perp \o_t-R^{-1}u_t\right|^2dx\\ &=\int\left|\n u_{tt}+\n_t\dot u+\n u\cdot \na u_t +\n u_t\cdot\na u-\na (\lm_t\div u)+\na P_t\right|^2dx\\&\le C\int\n|u_{tt}|^2dx+C\|\bar x^{(a+1)/2}\n_t\|_{L^q}^2\|\bar x^{-1}\dot u\|_{L^{2q/(q-2)}}^2 \\&\quad+C\|\na u_t\|_{L^2}^2+C\|\bar x^{-a} u_t\|_{L^4}^2\|\na u \|_{L^4}^2+C\|\na\lm_t\|_{L^2}^2\|\na u\|_{L^\infty}^2\\&\quad+C\|\lm_t\|_{L^q}^2\|\na^2u\|_{L^{2q/(q-2)}}^2+C
\|\na P_t\|_{L^2}^2\\&\le C\int\n|u_{tt}|^2dx +C\|\na u_t\|_{L^2}^4 +C,\ea\enn where in the last inequality we have used   \eqref{oo7}, \eqref{5.22}, and \eqref{o6}--\eqref{o7}. This combined with \eqref{lp} and  \eqref{o1} yields that
   \be \la{oo10}\ba  \|\na^2  u_t\|_{L^2} &\le C\|\na^\perp\o_t\|_{L^2}  +C\|\na\div   u_t\|_{L^2}\\ &\le C\|\na^\perp\o_t\|_{L^2} +C\|\na((2\mu+\lambda)\div
    u_t)\|_{L^2}\\&\quad+C \|\div   u_t\|_{L^{2q/(q-2)}} \|\na \lm\|_{L^q}   \\&\le C\|\n^{1/2} u_{tt}\|_{L^2}  +C\|\na u_t\|^2_{L^2} +C +\frac12  \|\na^2   u_t\|_{L^2}.\ea\ee
Putting \eqref{oo10} into \eqref{oo11} and choosing $\de$ suitably small lead to
 \be\la{oo9}\ba 2\Psi'(t)+\int\n|u_{tt}|^2dx\le  C \|\na u_t\|_{L^2}^4+C .\ea\ee

  Multiplying \eqref{oo9} by $t $   and
integrating it over $(0,T_0),$ we obtain from Gronwall's inequality, \eqref{o3}, and \eqref{oo1} that
\bnn \sup\limits_{0\le t\le T_0}t\left(\|\na u_t\|_{L^{2}}^2 +R^{-1}\|u_t\|_{L^{2}}^2\right) +\int_{0}^{T_0}t\|\n^{1/2}  u_{tt}\|_{L^{2}}^2 dt\le C , \enn which together with \eqref{oo10} and \eqref{oo1}  gives
\eqref{oo8} and finishes the proof of Lemma \ref{oiq1}.

\begin{lemma}\label{lem4.a5}  It holds that
\begin{equation}\la{5.13n}\ba
 \sup_{0\leq t\leq T_0}\left(\|\nabla^2 \n\|_{L^q }+\|\nabla^2  \lm \|_{L^q }+\|\nabla^2 P  \|_{L^q }\right) \leq C .\ea
\end{equation}
\end{lemma}

\pf Applying the differential operator $\na^2$ to both sides of
(\ref{a4.r3}), multiplying the resulting equations by $q|\na^2
g(\n)|^{q-2}\na^2 g(\n)$, and integrating it by parts over
$\O $ lead to
\be\la{4.r4}\ba
\left(\|\nabla^2 g \|_{L^q}\right)_t &\leq  C \left(\|\nabla
u\|_{L^\infty}\|\nabla^2 g \|_{L^q}+\|\nabla g \|_{L^\infty}\|\nabla^2
u\|_{L^q}+\|\nabla^2u\|_{W^{1,q}}\right)  \\
&\leq C \left(1+\|\nabla^2
u\|_{L^q}\right)
\left(1+\|\nabla^2g\|_{L^q} \right)+C \|\nabla^3u\|_{L^q} .
\ea\ee
By \eqref{lp}, the last term on the right-hand side of \eqref{4.r4} can be estimated as follows:
\be\label{4.f1}\ba
\|\nabla^3 u\|_{L^q}  &\le  C\left(\|\na(\na^\perp \o)\|_{L^q}+ \|\na(\na \div u)\|_{L^q} \right)+C\|\na^2 u\|_{L^q}\\ &\le  C\left(\| \na^2 \o \|_{L^q}+\|\na^2 F  \|_{L^q}+\|\na^2 P \|_{L^q}  \right)\\&\quad+C\left(\||\na u||\na^2 \lm| \|_{L^q}+\||\na^2u||\na \lm |\|_{L^q}\right)+C+\frac14\|\na^3u\|_{L^q}\\ &\le  C \|\na (\n\dot u)\|_{L^q}+C\|\na^2 P \|_{L^q}+C\|\na^2 \lm \|_{L^q}\|\na u\|_{L^\infty} +C\\&\quad+ \frac12\|\na^3u\|_{L^q}  ,
\ea\ee  where in the last inequality we have used \eqref{2.9-1} and the following simple fact:$$\||\na^2u||\na\lm|\|_{L^q}\le C\|\na^2u\|_{L^\infty}\|\na\lm\|_{L^q}\le  \ve\|\na^3u\|_{L^q}+C(\ve)$$ due to \eqref{oo2}. For the first term on the right-hand side of \eqref{4.f1}, it follows from the Sobolev inequality,  \eqref{o1},     \eqref{5.22},    \eqref{oo2},  and \eqref{5.b3} that
 \be\la{5.b4}\ba \|\na (\n\dot u)\|_{L^q}&\le C\|\bar x^{-a}\na \dot u \|_{L^q} +C\|\bar x^{-a}\dot u \|_{L^\infty}\|\bar x^a\na \n\|_{L^q}\\&\le C\|\bar x^{-a}\na \dot u \|_{L^q} +C\|\bar x^{-a}\dot u \|_{L^q}+ C\|\na(\bar x^{-a} \dot u )\|_{L^q}\\&\le C\|\bar x^{-a}\na \dot u \|_{L^q} +C\|\bar x^{-a}\dot u \|_{L^q}\\ &\le C\| \na  u_t \|_{L^q}+ C\|\bar x^{-a}|u||\na^2  u| \|_{L^q}+C\|\bar x^{-a}|\na  u|^2 \|_{L^q} \\&\quad +C\| \na  u_t \|_{L^2}+C \\ &\le C\| \na  u_t \|_{L^q} +C\| \na^3  u  \|_{L^2}+C\|\na u\|_{H^1}^2 +C\| \na  u_t \|_{L^2}+C  \\ &\le C\| \na  u_t \|^{2/q}_{L^2}\| \na^2  u_t \|^{1-2/q}_{L^2}  +C\| \na u_t  \|_{L^2}+C  ,\ea\ee
which together with \eqref{oo8} and  \eqref{oo1} yields
 \be\la{4.29}\ba &\int_0^{T_0} \|\na (\n\dot u)\|_{L^q}^{1+1/q}dt\\&\le C\int_0^{T_0}\left( (t\|\na   u_t\|_{L^2}^2)^{1/q}(t\|\na^2  u_t\|_{L^2}^2)^{(q-2)/(2q) }t^{-1/2}\right)^{1+1/q}dt+C\\ & \le C  \int_0^{T_0}\left(   t\|\na^2   u_t\|_{L^2}^2 + t^{-(q^2+q)/(q^2+q+2)}\right)dt+C \\ &\le  C. \ea\ee

  Putting  \eqref{4.f1} into \eqref{4.r4}, we obtain \eqref{5.13n} from Gronwall's inequality, \eqref{4.29}, and \eqref{o1}. The proof of Lemma \ref{lem4.a5} is completed.

\begin{lemma} \la{4.m44} It holds that
\begin{equation}\label{4.m30}\ba
&\sup_{0\leq t\leq T_0} t\left(  \|\na^3
u \|_{L^2\cap L^q} + \|\na
u_t \|_{H^1} +\|\na^2(\n u)\|_{L^{(q+2)/2}} \right)\\& +\int_0^{T_0}  t^2\left( \|\nabla u_{tt}\|_{L^2}^2 +\|u_{tt}\bar x^{-1}\|_{L^2}^2 +R^{-1}\|u_{tt}\|_{L^2}^2\right)dt\leq
C .\ea
\end{equation}

\end{lemma}

\pf  We claim that
\begin{equation}\la{5.o2}
\sup_{0\leq t\leq T_0}t^2 \|\rho^{1/2}
u_{tt}\|_{L^2}^2+\int_0^{T_0}  t^2\left( \|\nabla u_{tt}\|_{L^2}^2 +R^{-1}\|u_{tt}\|_{L^2}^2\right)dt  \leq
C ,
\end{equation}
which together with  \eqref{oo10},  \eqref{oo8},   and \eqref{3.ii2} yields  that
\begin{equation}\la{5.p6}
\sup_{0\leq t\leq T_0}t \|\na   u_t\|_{H^1} +\int_0^{T_0}t^2\|u_{tt}\bar x^{-1}\|_{L^2}^2dt\leq
C .
\end{equation}
This combined with  \eqref{5.b3},
\eqref{4.f1}, \eqref{5.b4}, and  \eqref{5.13n}  leads to
\begin{equation}\la{5.p5}
\sup_{0\leq t\leq T_0}t\|\na^3  u\|_{L^2\cap L^q}\leq
C ,
\end{equation}  which, along with    \eqref{o1},  \eqref{5.13n},   and \eqref{5.13a}, gives \be\la{5.p7}\ba & t\|\na^2(\n u)\|_{L^{(q+2)/2}}\\&\le Ct\||\na^2 \n|| u |\|_{L^{(q+2)/2}} +Ct\||\na  \n||\na u |\|_{L^{(q+2)/2}}+Ct\|\n\na^2  u \|_{L^{(q+2)/2}} \\& \le Ct\|\bar x^{\de_0}\na^2\n \|_{L^2}^{2/(q+2)}\|\na^2\n\|_{L^q}^{q/(q+2)}\|\bar x^{-2\de_0/(q+2)} u\|_{L^\infty} \\&\quad+Ct\|\na \n\|_{L^q}\|\na u\|_{L^{q(q+2)/(q-2)}} +Ct \|\na^2  u\|_{L^{(q+2)/2}} \le C.\ea\ee We thus directly obtain   \eqref{4.m30} from   \eqref{5.o2}-\eqref{5.p7}.

It remains to prove \eqref{5.o2}. In fact, differentiating   (\ref{zb1})   with
respect to $t$ leads to
\begin{eqnarray*}
&&\rho u_{ttt}+\rho u\cdot\nabla u_{tt}-\nabla((2\mu+\lambda)\divg u_{tt})-\mu\na^\perp \o_{tt}+R^{-1}u_{tt}\\
&&\quad=2\na(\lm_t\div u_t)+\na(\lm_{tt}\div u)+2\divg (\rho u)u_{tt}+\divg(\rho u)_tu_t \\
&&\qquad-2(\rho
u)_t\cdot\nabla u_t- \rho_{tt}u\cdot\nabla u-2\rho_tu_t \cdot\nabla u-\rho u_{tt}\cdot\nabla
u-\nabla P_{tt} ,
\end{eqnarray*}
which, multiplied by $u_{tt}$  and integrated by parts over
$\O $, yields
\be\ba
&\frac{1}{2}\frac{d}{dt}\int \rho|u_{tt}|^2dx+\int \left((2\mu+\lambda)(\div  u_{tt})^2+\mu  \o_{tt}^2+R^{-1}|u_{tt}|^2\right)dx \\
&=-2\int \lm_t\div u_t\div u_{tt}dx- \int \lm_{tt}\div u \div u_{tt}dx \\& \quad-4\int \rho u\cdot\nabla u_{tt}\cdot u_{tt}dx-\int (\rho
u)_t\cdot\left(\nabla(u_t\cdot
u_{tt})+2\nabla u_t\cdot u_{tt}\right)dx \\
&\quad -\int  (\rho u)_{t}\cdot \na (u  \cdot \nabla u\cdot
u_{tt})dx-2\int   \rho_tu_t \cdot \nabla u\cdot
u_{tt}dx\\
& \quad-\int \rho u_{tt}\cdot\nabla u\cdot
u_{tt}dx+\int P_{tt}\divg u_{tt}dx \triangleq\sum_{i=1}^8J_i.\label{4.31}
\ea\ee

We estimate each $J_i(i=1,\cdots,8)$   as
follows:

First,  we deduce from \eqref{oo7}    that
\be\la{j1}\ba |J_1|&\le  C  \|\lm_t\|_{L^q}\|\na u_t\|_{L^{2q/(q-2)}}\|\na u_{tt}\|_{L^2}\\&\le   \ve \|\na u_{tt}\|_{L^2}^2+ C(\ve) \|\na u_t\|_{H^1}^2 . \ea\ee

Next, the Cauchy inequality gives
\be\la{pl1} \ba|J_2|& \le   \ve \|\na u_{tt}\|_{L^2}^2 +C(\ve) \|\lm _{tt}\|_{L^2}^2\|\na u \|_{L^\infty}^2 \\& \le  \ve \|\na u_{tt}\|_{L^2}^2 +C(\ve) +C(\ve)\|\na u_t\|_{L^2}^4 +C(\ve) \|\lm _{tt}\|_{L^2}^4,\ea\ee where in the second inequality we have used   \eqref{5.b3}. Using \eqref{bb1}, we estimate the last term on the right-hand side of \eqref{pl1} as follows:
\be\la{5.a9}\ba \|\lm _{tt} \|_{L^2} & \le C\| |u_t||\na \lm | \|_{L^2}+C\| |u||\na \lm_t | \|_{L^2} +C\|  \lm_t \div u  \|_{L^2}  +C\|  \lm   \div u_t  \|_{L^2}\\& \le C\|\bar x^{-\te_0}u_t\|_{L^{2q/((q-2)\te_0)}} \|\bar x^{\te_0}\na \lm  \|_{L^{2q/(q-(q-2)\te_0)}}\\&\quad+C\|\bar x^{- \de_0/2}u\|_{L^\infty}\|\bar x^{  \de_0/2 }\na \lm_t \|_{L^2}  +C\|  \lm_t \|_{L^q}\|\div u \|_{L^{2q/(q-2)}}  \\&\quad+C\|    \na u_t  \|_{L^2}\\&\le C+C\|    \na u_t  \|_{L^2},  \ea\ee where in the last inequality we have used \eqref{oo2}-\eqref{5.22}, \eqref{o1}, and the following simple fact that
\bnn\ba  \|\bar x^{  \de_0/2 }\na \lm_t \|_{L^2} &\le C\| \bar x^{ \de_0/2}|u| |\na^2\lm| \|_{L^2}+C\|\bar x^{\te_0}|\na u||\na \lm|\|_{L^2}+C\|\bar x^{\beta a}\lm \na^2 u\|_{L^2}\\ &\le C\| \bar x^{ -\de_0/2} u\|_{L^\infty} \|\bar x^{\de_0}\na^2\lm  \|_{L^2}+C\| \na u\|_{L^{2q/(q-2)}}\|\bar x^{\te_0}\na \lm \|_{L^q}\\&\quad +C\|\bar x^{\beta a}\lm\|_{L^\infty}\| \na^2 u\|_{L^2}\\ &\le C, \ea\enn
due to \eqref{oo3},  \eqref{5.13a},  \eqref{oo2}, and  \eqref{o1}.
Putting  \eqref{5.a9} into  \eqref{pl1}  gives
\be \la{j2}\ba|J_2|  \le   \ve \|\na u_{tt}\|_{L^2}^2 +C(\ve) (1+\|\na u_t\|_{L^2}^4) .\ea\ee

Next, the combination of the Cauchy inequality with \eqref{oo3}  yields that
\be \la{j3}\ba
|J_3|&\leq C\|\n^{1/2}u\|_{L^{\infty}}\|\rho^{1/2} u_{tt}\|_{L^2}\|\nabla u_{tt}\|_{L^2}\leq
\ve \|\nabla
u_{tt}\|_{L^2}^2+C(\ve)\|\rho^{1/2} u_{tt}\|_{L^2}^2.\ea\ee

Next, noticing that   \eqref{oo3}-\eqref{5.22} lead to
\be \la{oo32}\ba  \|\bar x  (\n u)_t \|_{L^q}   &\le C\|\bar x |\n_t||u|\|_{L^q}+C\|\bar x \n|u_t|\|_{L^q}\\&\le C \|\n_t\bar x^{(a+1)/2}\|_{L^q} \|\bar x^{(1-a)/2}u\|_{L^\infty}
\\&\quad +C\|\n\bar x^a\|_{L^{2q/(3-a)}}\|u_t\bar x^{1-a}\|_{L^{2q/(a-1)}} \\&\le C +C\|\na u_t\|_{L^2} ,\ea\ee we obtain from Holder's   inequality,   \eqref{oo7},   \eqref{3.i2},   and \eqref{3.v2}  that
\be\la{j4}\ba |J_4|&\le C\|\bar x  (\n u)_t \|_{L^q}  \left(\|\bar x^{-1}u_{tt}\|_{L^{2q/(q-2)}}\|\na u_t\|_{L^2}+\|\bar x^{-1}u_t\|_{L^{2q/(q-2)}}\|\na u_{tt}\|_{L^2}\right)\\ &\le C \left(1+\|\na u_t\|_{L^2}^2\right)\left(\|\n^{1/2}u_{tt}\|_{L^2}+\|\na u_{tt}\|_{L^2}\right)\\ &\le C(\ve) \left(1+\|\na u_t\|_{L^2}^4\right)+\ve\left(\|\n^{1/2}u_{tt}\|^2_{L^2}+\|\na u_{tt}\|^2_{L^2}\right) .\ea\ee

Then, it follows from   \eqref{oo32},   \eqref{oo3},  \eqref{oo2},    and \eqref{3.v2}    that
\be\la{j5}\ba |J_5|&\le C\int |(\n u)_t| ( |u ||\na^2 u|| u_{tt}|+|u||\na u||\na u_{tt}|+|\na u|^2|u_{tt}|)dx\\ &\le C \|\bar x  (\n u)_t \|_{L^q}\|\bar x^{-1/q}u\|_{L^\infty}\|\na^2u\|_{L^2}\|\bar x^{-(q-1)/q}u_{tt}\|_{L^{2q/(q-2)}}\\&\quad+ C \|\bar x  (\n u)_t \|_{L^q}\|\bar x^{-1}u\|_{L^\infty}\|\na u\|_{L^{2q/(q-2)}}\|\na u_{tt}\|_{L^2}\\&\quad + C \|\bar x  (\n u)_t \|_{L^q}\|\na u \|_{L^4}^2\|\bar x^{-1} u_{tt}\|_{L^{2q/(q-2)}}\\ &\le C \left(1+\|\na u_t\|_{L^2}\right)\left(\|\n^{1/2}u_{tt}\|_{L^2}+\|\na u_{tt}\|_{L^2}\right)\\ &\le C(\ve) \left(1+\|\na u_t\|_{L^2}^2\right)+\ve\left(\|\n^{1/2}u_{tt}\|^2_{L^2}+\|\na u_{tt}\|^2_{L^2}\right).\ea\ee

Next, Cauchy's inequality together with   \eqref{oo7},   \eqref{5.22},   and \eqref{3.v2} gives
\be\la{j6}\ba |J_6|&\le C \int |\n_t||u_t||\na u||u_{tt}| dx\\& \le C\|\bar x\n_t\|_{L^q}\|\bar x^{-1/2}u_t\|_{L^{4q/(q-2)}}\|\na u\|_{L^2}\|\bar x^{-1/2}u_{tt}\|_{L^{4q/(q-2)}}\\ &\le C \left(1+\|\na u_t\|_{L^2} \right)\left(\|\n^{1/2}u_{tt}\|_{L^2}+\|\na u_{tt}\|_{L^2}\right)\\ &\le C(\ve) \left(1+\|\na u_t\|_{L^2}^2\right)+\ve\left(\|\n^{1/2}u_{tt}\|^2_{L^2}+\|\na u_{tt}\|^2_{L^2}\right) .\ea\ee

Finally,  similar to  \eqref{5.a9}, we have
\bnn \ba   \| P_{tt} \|_{L^2}    \le  C(1+  \|  \na u_t \|_{L^2}) , \ea\enn which together with direct calculations gives
\be\la{j7}\ba  |J_7|+|J_8|&\le C\|\na u\|_{L^\infty}\int \n |u_{tt}|^2 dx+\ve\int| \na u_{tt}|^2 dx +C(\ve)\|P_{tt}\|_{L^2}^2\\ &\le C\|\na u\|_{L^\infty}\|\n^{1/2}u_{tt}\|^2_{L^2}+\ve\|\na u_{tt}\|^2_{L^2}
  +C(\ve)(1 +\|  \na  u_t \|_{L^2}^2) . \ea\ee

Substituting \eqref{j1},  \eqref{j2},  \eqref{j3},  and  \eqref{j4}-\eqref{j7}   into \eqref{4.31}, choosing $\ve$ suitably small, and multiplying the resulting inequality by $t^2,$ we obtain \eqref{5.o2} after  using Gronwall's inequality and \eqref{oo8}.
 The proof of Lemma \ref{4.m44} is  finished.

\section{Proofs of Theorems \ref{t1} and \ref{t2} }

To prove Theorems \ref{t1} and \ref{t2},
we will only deal with the case that $\beta>0,$ since  the same  procedure can be applied to  the case that $\beta=0$ after some small modifications.

{\it Proof of Theorem  \ref{t1}.}
Let  $(\n_{0},u_{0})$  be as  in Theorem \ref{t1}. Without loss of generality, assume that
\bnn \int_{\rr} \n_0dx=1,\enn  which implies that there exists a positive constant $N_0$ such that  \be\la{oi3.8} \int_{B_{N_0}}  \n_0  dx\ge \frac34\int_{\rr}\n_0dx=\frac34.\ee
We
construct
$\n_{0}^{R}=\hat\n_{0}^{R}+R^{-1}e^{-|x|^2} $ where  $0\le\hat\n_{0}^{R}\in  C^\infty_0(\rr)
$  satisfies that \be \la{bbi1} \int_{B_{N_0}}\hat\n^R_0dx\ge 1/2,\ee
and that \be\la{bci0}\begin{cases} \bar x^a \hat\n_{0}^{R}\rightarrow \bar x^a \n_{0}\quad {\rm in}\,\, L^1(\rr)\cap H^{1}(\rr)\cap W^{1,q}(\rr) ,\\  \lm(\hat\n_{0}^{R})\rightarrow  \lm( \n_{0}),\quad \bar x^{\te_0 }\na\lm(\hat\n_{0}^{R})\rightarrow \bar x^{\te_0 }\na\lm( \n_{0})\quad {\rm in}\,\, L^2(\rr)\cap L^q(\rr) , \end{cases}\ee as $R\rightarrow\infty.$

Since $\na u_0\in L^2(\rr),$  choosing $v^R_i\in C^\infty_0(\O)(i=1,2 )$ such that \be \la{bci3}\lim\limits_{R\rightarrow \infty}\|v^R_i-\pa_iu_0\|_{L^2(\rr)}=0,\quad i=1,2,\ee we consider the unique smooth solution $u_0^R$ of the following elliptic problem:
 \be \la{bbi2} \begin{cases}-  \lap u_{0}^{R} + R^{-1} u_0^R  =-\n_0^R u_0^R+\sqrt{\n_{0}^{R}} h^R-  \p_iv^R_i ,& {\rm in} \,\,  B_{R},\\ u_{0}^{R}\cdot n=0,\,\,{\rm rot}u^R_0=0,&{\rm on} \,\,\partial B_{R},\end{cases} \ee
where $  h^R= (\sqrt{\n_0}u_0  )*j_{1/R}   $  with $j_\de$
being the standard mollifying kernel of width $\de.$ Extending  $u_{0}^{R} $ to $\rr$ by defining $0$ outside $B_{R}$ and denoting $w_0^R\triangleq u_0^R\vp_R$ with $\vp_R$ as in \eqref{vp1}, we claim that
 \be \la{3.7i4} \lim\limits_{R\rightarrow \infty}\left(\|\na ( w_0^{R }-u_0)\|_{L^2(\rr)}+\|\sqrt{\n_0^R}  w_0^{R }-\sqrt{\n_0}u_0 \|_{L^2(\rr)}\right)=0.\ee

In fact, multiplying \eqref{bbi2} by $u_0^R$ and integrating the resulting equation over $\O$ lead to \bnn  \ba&\int_{\O} \left(\n_0^R+R^{-1}\right)|u_0^R|^2dx+ \int_{\O} (\rot u_0^R)^2dx +\int_{\O} (\div u_0^R)^2dx\\&\le \|\sqrt{\n_0^R} u_0^R\|_{L^2(\O)} \|h^R\|_{L^2(\O)} +C\|v_i^R\|_{L^2}\|\p_i u^R_0\|_{L^2}\\&\le \ve \|\na u_0^R\|_{L^2(\O)}^2+\ve \int_{\O} \n_0^R |u_0^R|^2dx+C(\ve) ,\ea\enn which implies \be \la{2.i9-4}  R^{-1}\int_{\O}|u_0^R|^2dx+\int_{\O} \n_0^R|u_0^R|^2dx+ \int_{\O} |\na u_0^R|^2dx  \le C,\ee for some $C$ independent of $R.$

 We deduce from   \eqref{2.i9-4} and   \eqref{bci0}  that there exists a subsequence $R_j\rightarrow \infty$ and a function $w_0\in\{w_0\in H^1_{\rm loc}(\rr)|\sqrt{\n_0}w_0\in L^2(\rr),   \na w_0\in L^2(\rr)\}$ such that
 \be\la{bci9}\begin{cases}\sqrt{\n^{R_j}_0}w^{R_j}_0 \rightharpoonup \sqrt{\n_0} w_0 \mbox{ weakly in } L^2(\rr) ,\\
\na w_0^{R_j}\rightharpoonup \na w_0 \mbox{ weakly in } L^2(\rr).\end{cases}\ee It follows from    \eqref{bbi2}    and \eqref{2.i9-4}   that $w_{0}^{R}$ satisfies
 \be \la{bci2} \ba  -  \lap w_{0}^{R}  + R^{-1} w_0^R   =-\n_0^R w_0^R+\sqrt{\n_{0}^{R}} h^R\vp_R -\p_iv^R_i\vp_R+R^{-1}F^R,\ea \ee
 with $\|F^R\|_{L^2(\rr)}\le C.$ Thus, one can deduce    from \eqref{bci2},  \eqref{bci9}, and \eqref{bci3}  that,  for any $\psi\in C_0^\infty(\rr),$
 \bnn \int_{\rr}\p_i(w_0-u_0)\cdot\p_i\psi dx+\int_{\rr} \n_0(w_0-u_0)\cdot\psi dx=0,\enn
 which   yields that \be\la{bai1} w_0=u_0.\ee Furthermore, we get from \eqref{bci2} that \bnn\ba  \limsup\limits_{R_j\rightarrow \infty}\int_{\rr}\left( |\na w_0^{R_j}|^2  + \n_0^{R_j}|w_0^{R_j}|^2\right)dx \le \int_{\rr}\left( |\na u_0 |^2  + \n_0 |u_0 |^2\right)dx,\ea\enn which combined with \eqref{bci9}  implies \bnn \lim\limits_{R_j\rightarrow \infty}\int_{\rr}|\na w_0^{R_j} |^2dx=\int_{\rr}|\na u_0 |^2dx,\,\,\lim\limits_{R_j\rightarrow \infty}\int_{\rr} \n_0^{R_j}   |w_0^{R_j} |^2dx=\int_{\rr}\n_0 | u_0 |^2dx.\enn
 This, along with  \eqref{bai1} and \eqref{bci9}, gives \eqref{3.7i4}.

Then, in terms of Lemma \ref{th0}, the initial-boundary-value problem  \eqref{b2} with the initial data $(\n_0^R,u_0^R)$ has  a classical solution   $(\n^{R},u^{R}) $ on $B_{R}\times [0,T_R].$  Moreover,  Proposition \ref{pro} shows that there exists a $T_0$ independent of $R$ such that \eqref{o1} holds for $(\n^{R},u^{R}) .$
Extending   $( \n^{R},u^{R})$  by zero on $\rr\setminus B_{R}$  and denoting $$\nr\triangleq (\vp_R)^{4/\ti\beta}\n^R,\quad w^R\triangleq \vp_Ru^R,$$ with $\vp_R$ as in \eqref{vp1} and $\ti\beta=\min\{\beta,1\},$ we first deduce from  \eqref{o1} that
\be\la{kq1} \ba  &\sup\limits_{0\le t\le T_0}\left(\|\sqrt{\ti\n^R } w^R\|_{L^2(\rr)}+\|\na w^R \|_{L^2(\rr)}\right)\\&\le  C+ C\sup\limits_{0\le t\le T_0}\left(\| \na u^R \|_{L^2(\O)}+C R^{-1}\|  u^R \|_{L^2(\O)}\right)\\&\le C ,\ea \ee
and that \be \ba&\sup\limits_{0\le t\le T_0}\left(\|\ti\n^R\bar x^a\|_{L^1(\rr)\cap L^\infty(\rr)}+\|\lm(\ti\n^R )\|_{L^2(\rr)}\right) \\&\le \sup\limits_{0\le t\le T_0}\left(\|\n^R\bar x^a\|_{L^1(\O)\cap L^\infty(\O)}+\|\lm(\n^R) \|_{L^2(\O)}\right) \le C.\ea\ee

Next, for $p\in [2,q] ,$ it follows from \eqref{b2a} and \eqref{o1} that
\be \ba & \sup\limits_{0\le t\le T_0}\left(\|\na (\bar x^{\te_0}  \lm(\ti\n^R) ) \|_{L^p(\rr)} + \|  \bar x^{\te_0} \na  \lm(\ti\n^R  ) \|_{L^p(\rr)}  \right)\\ &\le   C\sup\limits_{0\le t\le T_0}\left(\|\na(\bar x^{\te_0} \lm( \n^R)) \|_{L^p(\O)}+ \|(\bar x^a \n^R)^\beta \na\vp_R \|_{L^p(\O)}\right)\\&\quad +C\sup\limits_{0\le t\le T_0} \|  \bar x^{\te_0} \na  \lm( \n^R  ) \|_{L^p(\O)} \\& \le C  +C\sup\limits_{0\le t\le T_0}\|\bar x^a \n^R\|_{L^{2p\beta/\ti\beta}(\O)}^\beta  \|\na\vp_R\|_{L^{2p /(2-\ti\beta)}(\O)}\\&\le C, \ea \ee and that
  \be \ba &\sup\limits_{0\le t\le T_0}\left(\|\na(\bar x^a\ti \n^R)\|_{L^p(\rr)}+\|\bar x^a\na \ti \n^R \|_{L^p(\rr)}\right)\\&\le   C\sup\limits_{0\le t\le T_0}\left(\|\bar x^a\na \n^R \|_{L^p(\O)}+   \|\bar x^a \n^R \na  \vp_R  \|_{L^p(\O)}+  \|\n^R \na \bar x^a   \|_{L^p(\O)}\right)\\& \le   C +C \sup\limits_{0\le t\le T_0}\|\bar x^a \n^R \|_{L^p(\O)}\le C. \ea \ee

Then, it follows from \eqref{o1} and \eqref{3.22} that
\be \int_0^{T_0} \left(\norm[L^{q}(\rr)]{\nabla^2w^R}^{(q+1)/q}+t\|\na^2w^R\|_{L^q(\rr)}^{2}+ \|\na^2w^R\|_{L^2(\rr)}^{2}\right) dt\le C,\ee
and that for $p\in [2,q],$
\be\ba \int_0^{T_0}\|\bar x\ti\n^R_t\|^2_{L^p(\rr)}dt&\le C\int_0^{T_0}\left(\|\bar x |u^R||\na \n^R| \|^2_{L^p(\O)}+\|\bar x  \n^R \div u^R \|^2_{L^p(\O)}\right)dt\\ &\le C\int_0^{T_0} \|\bar x^{1-a} u^R\|_{L^\infty(\O)}^2\|\bar x^a\na \n^R  \|^2_{L^p(\O)}dt+C\\&\le C.
\ea\ee
Next, one derives from \eqref{li-1a} and \eqref{o1} that
\be\la{kq2}\ba  &\sup\limits_{0\le t\le T_0}t\int_{\rr}\ti\n^R|w_t^R|^2dx+\int_0^{T_0}t\|\na w^R_t\|_{L^2(\rr)}^2dt\\&\le C+ C \int_0^{T_0}t\left(\|\na u^R_t\|_{L^2(\O)}^2 + R^{-2} \|  u^R_t\|_{L^2(\O)}^2\right)dt\\&\le C.\ea\ee

  With all these estimates \eqref{kq1}-\eqref{kq2} at hand,  we  find  that the sequence
$(\ti\n^R,w^R)$ converges, up to the extraction of subsequences, to some limit $(\n,u)$   in the
obvious weak sense, that is, as $R\rightarrow \infty,$ we have
\be \la{kq3} R^{-1}w^R \rightarrow  0,\mbox{   in }L^2(\rr\times(0,T_0) ),\ee  \be \bar x\ti\n^R\rightarrow \bar  x\n, \mbox{ in } C(\overline{ B_N}\times [0,T_0]), \mbox{ for any }N>0,\ee
\be \bar x^a\ti\n^R\rightharpoonup \bar x^a \n ,\mbox{ weakly * in }L^\infty(0,T_0; H^1(\rr)\cap W^{1,q}(\rr)),\ee
\be \na(\bar x^{\te_0}\lm(\ti\n^R))\rightharpoonup \na(\bar x^{\te_0}\lm( \n)) ,\mbox{ weakly * in }L^\infty(0,T_0; L^2(\rr)\cap L^{ q}(\rr)),\ee\be  \sqrt{\ti\n^R} w^R\rightharpoonup \sqrt{\n } u,\,\, \na w^R\rightharpoonup \na u ,\mbox{ weakly * in }L^\infty(0,T_0; L^2(\rr)),\ee \be  \na^2 w^R\rightharpoonup \na^2 u ,\mbox{ weakly  in }L^{(q+1)/q}(0,T_0; L^q(\rr))\cap L^2(\rr\times (0,T_0)),\ee \be  t^{1/2}\na^2 w^R\rightharpoonup t^{1/2}\na^2 u ,\mbox{ weakly  in }L^2(0,T_0; L^q(\rr)),\ee  \be  \sqrt{t } \sqrt{\ti \n^R} w^R_t\rightharpoonup  \sqrt{t }\sqrt{ \n } u_t,\,\, \na w^R\rightharpoonup \na u ,\mbox{ weakly * in }L^\infty(0,T_0; L^2(\rr)),\ee  \be \sqrt{ t} \na  w^R_t\rightharpoonup \sqrt{t }\na  u_t ,\mbox{ weakly  in } L^2(\rr\times (0,T_0)),\ee
with \be \la{kq4} \bar x^a \n \in L^\infty(0,T_0; L^1(\rr)), \quad \inf\limits_{0\le t\le T_0}\int_{B_{2N_0}}\n(x,t)dx\ge \frac14. \ee

Next, for any function $\phi\in C^\infty_0(\rr\times [0,T_0)),$  we take $\phi(\vp_R)^{8/\ti\beta}$ as test function in the initial-boundary-value problem  \eqref{b2} with the initial data $(\n_0^R,u_0^R).$ Then letting $R\rightarrow \infty,$ it follows from \eqref{kq3}-\eqref{kq4} that
$(\n,u)$ is a strong solution of \eqref{n1}-\eqref{n4} on $\rr\times (0,T_0]$ satisfying \eqref{1.10} and \eqref{l1.2}.  The  proof of the existence part of Theorem \ref{t1} is finished.

It only remains to prove the uniqueness of the strong solutions  satisfying \eqref{1.10} and \eqref{l1.2}. We only treat the case  $\beta>0,$ since the   procedure  can be adapted  to  the case $\beta=0$ after some small modifications.   Let   $(\rho,u)$ and $(\bar\rho,\bar u)$ be two strong solutions satisfying \eqref{1.10} and \eqref{l1.2}  with the same initial data.
   Subtracting the momentum equations satisfied by $(\rho,u)$ and $(\bar\rho,\bar u)$ yields
\be\ba\la{5.5}
& \rho U_t + \rho u\cdot\nabla U -\mu\lap U - \na\left((\mu+\lambda (\n))\div  U\right)\\
 &=  - \rho U\cdot\nabla\bar u-H(\bar u_t+\bar u\cdot\nabla\bar u)-\na \left( P(\n)-P(\bar{\n}) \right) +\na\left((\lambda (\n)-\lambda(\bar{\n}))\div \bar{u}\right),
\ea\ee with
$$ H\triangleq\rho-\bar\rho, \quad U\triangleq u-\bar u.$$
Since $\mu+\lm\ge 0,$ multiplying (\ref{5.5}) by $U$ and integrating by parts lead to
\be\ba\la{5.6}
& \frac{d}{dt}\int \rho |U|^2dx  + 2\mu \int |\na U|^2 dx\\
 & \le   C \norm[L^{\infty}]{\nabla \bar{u}}\int \n |U|^2dx+C\int|H || U| \left(| \bar u_{t}|+ | \bar u ||\na\bar{u}|\right) dx  \\
 &\quad+C\left(\norm[L^{2}]{ P(\n)-P(\bar{\n})} +\norm[L^{\infty}]{\na \bar{u}}\norm[L^{2}]{\lambda (\n)-\lambda(\bar{\n})}\right)\norm[L^{2}]{\div U}\\
 & \triangleq   C \norm[L^{\infty}]{\nabla \bar{u}}\int \n |U|^2dx+K_1+K_2.
\ea\ee

We first estimate $K_1.$   Holder's inequality shows that for $r\in (1,a),$
  \be\ba\la{v5.8}
   K_1& \le C\norm[L^2]{H \bar{x}^{r}}\norm[L^{4}] {U \bar{x}^{-r/2}}\left( \norm[L^{4}]{\bar u_{t}\bar{x}^{-r/2}}+\norm[L^{\infty}]{\na \bar{u}} \norm[L^{4}] {\bar u\bar{x}^{-r/2}}\right) \\
&\le C(\ve)\left( \norm[L^{2}] {\sqrt{\bar{\n}}\bar{u}_{t}}^{2}+ \norm[L^{2}] {\na\bar{u}_{t}}^{2}+\norm[L^{\infty}]{\na \bar{u}}^{2}\right)\norm[L^2]{H \bar{x}^{r}}^2\\ & \quad+\ve\left(\norm[L^{2}]{\sqrt{\rho} U}^{2}+ \norm[L^{2}]{\na U}^2\right) ,
 \ea\ee
 where in the second inequality we have used  \eqref{3.i2} and   \eqref{l1.2}. Then, subtracting the mass   equation  for $(\rho,u)$ and $(\bar\rho,\bar u)$ gives
\be\la{5.2}
 H_t + \bar u\cdot\nabla H +H\div \bar u + \rho \div U  + U\cdot\nabla \rho= 0.
\ee
 Multiplying   \eqref{5.2} by $2H\bar{x}^{2r}$ and  integrating   by parts lead to
\bnn
\ba
 &\left(\norm[L^{2}]{H\bar{x}^{r}}^{2}\right)_t\\
 &\le C\left(\norm[L^{\infty}]{\nabla \bar{u}}+\norm[L^{\infty}]{ \bar{u} \bar{x}^{-1/2}}\right)\norm[L^{2}]{H\bar{x}^{r}}^{2}+  C \norm[L^{\infty}]{\rho\bar{x}^{r}}\norm[L^{2}]{\nabla U}\norm[L^{2}]{H\bar{x}^{r}}\\
 &\quad+C\norm[L^{2}]{H\bar{x}^{r}}\norm[L^{2q/((q-2)(a-r))}]{ U  \bar{x}^{-(a-r)}}\norm[L^{2q/(q-(q-2)(a-r))}]
 { \bar{x}^a \na \n} \\
 &\le   C\left( 1+\norm[ W^{1,q}]{\na\bar u}  \right)\norm[L^{2}]{H\bar{x}^{r}}^{2}+C\norm[L^{2}]{H\bar{x}^r}
\left(\norm[L^{2}]{\nabla U}+\norm[L^{2}]{\sqrt{\n} U}\right),
\ea
\enn
 where in the second inequality we have used \eqref{l1.2}, \eqref{3.v2}, and \eqref{3.22}. This combined with
 Gronwall's inequality yields that  for all $0\le t\le T_{0}$
\be\ba\la{5.1}
\norm[L^{2}]{H\bar{x}^{r}}
\le & C\int_{0}^{t}\left(\norm[L^{2}]{\nabla U}+\norm[L^{2}]{\sqrt{\n} U}\right)ds.
\ea\ee

 As observed by Germain \cite{Ge}, putting \eqref{5.1} into \eqref{v5.8}
leads to
\be\ba\la{5.8}
   K_1& \le
C(\ve)\left(1+t\norm[L^{2}] {\nabla  {\bar{u}_t}}^{2}+t\norm[L^{q}]{ \na^{2} \bar{u}}^{2} \right)\int_{0}^{t} \left(\norm[L^{2}]{\nabla U}^{2}+\norm[L^{2}]{\sqrt{\n} U}^{2}\right)ds\\
&\quad+\ve \left(\norm[L^{2}]{\sqrt{\rho} U}^{2}+ \norm[L^{2}]{\na U}^2\right).
 \ea\ee

Next,  we will estimate $K_2.$ In fact, one deduces from \eqref{bb1} that
\bnn\ba  &(\lm(\n)-\lm(\bar\n))_t+
\bar u\cdot\na(\lm(\n)-\lm(\bar\n))\\&+U\cdot \na\lm( \n)+\beta(\lm(\n)-\lm(\bar\n))\div\bar u+\beta\lm( \n)\div U=0,\ea\enn
which gives \be\la{w2}\ba (\|\lm(\n)-\lm(\bar\n)\|_{L^2})_t\le & C(1+\|\na \bar u\|_{L^\infty})\|\lm(\n)-\lm(\bar\n)\|_{L^2}\\&+C\|U\cdot \na \lm(\n)\|_{L^2}+C\|\na U\|_{L^2}.\ea\ee
It follows from \eqref{1.10}, \eqref{l1.2}, and \eqref{3.i2} that
\bnn\ba \|U\cdot \na \lm(\n)\|_{L^2}&\le
 \|U\bar x^{-\te_0}\|_{L^{2q/((q-2)\te_0)}}\|\bar x^{\te_0}\na\lm(\n)\|_{L^{2q/(q-(q-2)\te_0)}}\\& \le C\norm[L^{2}]{\nabla U}+C\norm[L^{2}]{\sqrt{\n} U},\ea\enn which together with \eqref{w2} and Gronwall's inequality gives \be\la{w3}\norm[L^{2}]{\lambda (\n)-\lambda(\bar{\n})}\le C \int_{0}^{t}\left(\norm[L^{2}]{\nabla U}+\norm[L^{2}]{\sqrt{\n} U}\right)ds.\ee
Similarly,  we have
 \bnn \norm[L^{2}]{ P(\n)-P(\bar{\n})} \le C \int_{0}^{t}\left(\norm[L^{2}]{\nabla U}+\norm[L^{2}]{\sqrt{\n} U}\right)ds,\enn
which combined with \eqref{w3} shows   \be\ba\la{5.7-7}
 K_2\le \ve \norm[L^{2}]{\nabla U}^{2} +C(\ve)\left(1+ t\norm[L^{q}]{ \na^{2} \bar{u}}^{2}\right)  \int_{0}^{t}(\norm[L^{2}]{\nabla U}^{2}+\norm[L^{2}]{\sqrt{\rho} U}^{2})ds.
\ea\ee

Denoting $$G(t)\triangleq   \norm[L^{2}]{\sqrt{\rho}U}^{2}+  \int_{0}^{t} \left(\norm[L^{2}]{\sqrt{\rho}U}^{2}+\mu \|\na U\|_{L^2}^2\right)  ds,$$ putting \eqref{5.8} and \eqref{5.7-7} into  (\ref{5.6}) and choosing  $\ve$ suitably small  lead to
\bnn\ba
  G'(t) \le   C  \left(1+\norm[L^{\infty}]{\nabla \bar{u}}+t\norm[L^{q}]{ \na^{2} \bar{u}}^{2} + t\norm[L^{2}]
{\na \bar{u}_{t}}^{2} \right) G,
\ea\enn
 which together with  Gronwall's  inequality and \eqref{1.10} yields $G(t)=0.$ Hence, $U(x,t)=0 $ for  almost everywhere $(x,t)\in \rr\times(0,T_0).$  Then, \eqref{5.1} implies that $H(x,t)=0 $ for almost everywhere $(x,t)\in \rr\times(0,T_0).$ The proof of Theorem  \ref{t1} is completed.

{\it Proof of Theorem  \ref{t2}.}
Let  $(\n_{0},u_{0})$  be as  in Theorem \ref{t2}. Without loss of generality, assume that
\bnn \int_{\rr} \n_0dx=1,\enn  which implies that there exists a positive constant $N_0$ such that  \eqref{oi3.8} holds.
We construct
$\n_{0}^{R}=\hat\n_{0}^{R}+R^{-1}e^{-|x|^2} $ where  $0\le\hat\n_{0}^{R}\in  C^\infty_0(\rr)
$  satisfies \eqref{bbi1}, \eqref{bci0},  and  \be\la{bc10}\begin{cases} \na^2 \hat\n_{0}^{R} \rightarrow  \na^2   \n_{0} ,\, \na^2 \lm(\hat\n_{0}^{R})\rightarrow \na^2  \lm( \n_{0}),\, \na^2 P(\hat\n_{0}^{R})\rightarrow \na^2 P( \n_{0}),\quad {\rm in }\,  L^q(\rr)  ,\\ \bar x^{\de_0}\na^2 \hat\n_{0}^{R} \rightarrow  \bar x^{\de_0}\na^2   \n_{0} ,\,\, \bar x^{\de_0}\na^2 \lm(\hat\n_{0}^{R})\rightarrow \bar x^{\de_0}\na^2  \lm( \n_{0}),\quad {\rm in }\,  L^2(\rr) ,\\ \bar x^{\de_0}\na^2 P(\hat\n_{0}^{R})\rightarrow \bar x^{\de_0}\na^2 P( \n_{0})  ,\quad {\rm in }\, L^2(\rr) , \end{cases}\ee as $R\rightarrow\infty.$

Then, we consider the unique smooth solution $u_0^R$ of the following elliptic problem:
 \be \la{bb12} \begin{cases}-\mu \lap u_{0}^{R}-\na \left((\mu+\lambda(\n_{0}^{R}))\div u_{0}^{R}\right)+\na P(\n_{0}^{R})+ R^{-1} u_0^R \\ \quad =-\n_0^R u_0^R+\sqrt{\n_{0}^{R}} h^R,& {\rm in} \,\,  B_{R},\\ u_{0}^{R}\cdot n=0,\,\,{\rm rot}u^R_0=0,&{\rm on} \,\,\partial B_{R},\end{cases} \ee
where $  h^R= (\sqrt{\n_0}u_0+g )*j_{1/R}   $  with $j_\de$
being the standard mollifying kernel of width $\de.$
Multiplying \eqref{bb12} by $u_0^R$ and integrating the resulting equation over $\O$ lead to \bnn  \ba&\int_{\O} \left(\n_0^R+R^{-1}\right)|u_0^R|^2dx+\mu\int_{\O} |\rot u_0^R|^2dx +\int_{\O}(2\mu+\lambda(\n_{0}^{R}))(\div u_0^R)^2dx\\&\le \int_{\O}P(\n_0^R)|\div u_0^R|dx+\|\sqrt{\n_0^R} u_0^R\|_{L^2(\O)} \|h^R\|_{L^2(\O)}\\&\le \ve \|\na u_0^R\|_{L^2(\O)}^2+\ve \int_{\O}  \n_0^R |u_0^R|^2dx+C(\ve),\ea\enn which implies \be \la{2.9-4}  R^{-1}\int_{\O}|u_0^R|^2dx+\int_{\O} \n_0^R|u_0^R|^2dx+ \int_{\O} |\na u_0^R|^2dx  \le C,\ee for some $C$ independent of $R.$
By \eqref{lp}, we have  \bnn\ba    \|\na^2u_0^R\|_{L^2(\O)}
 &\le C  \|  \na {\rm rot}u_0^R\|_{L^2(\O)} +  C \|\na\left((2\mu+\lambda(\n_0^R))\div u_0^R\right)\|_{L^2(\O)}\\ &\quad +  C\||\na \lambda(\n_0^R)|\div u_0^R\|_{L^{2}(\O)} \\&\le C\|\mu \lap u_{0}^{R}+\na \left((\mu+\lambda(\n_{0}^{R}))\div u_{0}^{R}\right)- R^{-1} u_0^R \|_{L^2} \\&\quad  +C\|\na \lambda(\n_0^R)\|_{L^{q}(\O)}   (1+\|\na^2 u^R_0\|_{L^2(\O)}^{2/q} ) \\
&\le  C  \|\n_0^R  u_0^R\|_{L^2(\O)}  + C\|\na P(\n_0^R)\|_{L^2(\O)}\\&\quad +C\|\sqrt{\n_0^R}  h^R\|_{L^2(\O)} +C+\frac12\|\na^2u_0^R\|_{L^2(\O)}\\&\le C+\frac12\|\na^2u_0^R\|_{L^2(\O)},\ea\enn which gives \be \la{2.9-3}\|\na^2u_0^R\|_{L^2(\O)}\le C.\ee

   Next, extending  $u_{0}^{R} $ to $\rr$ by defining $0$ outside $B_{R}$ and denoting $w_0^R\triangleq u_0^R\vp_R$ with $\vp_R$ as in \eqref{vp1},
 we deduce from \eqref{2.9-4} and \eqref{2.9-3} that  \bnn \|\na w_0^R\|_{H^1(\rr)} \le C ,\enn
which together with \eqref{2.9-4} and  \eqref{bc10}  yields that there exists a subsequence $R_j\rightarrow \infty$ and a function $w_0\in\{w_0\in H^2_{\rm loc}(\rr)|\sqrt{\n_0}w_0\in L^2(\rr),   \na w_0\in H^1(\rr)\}$ such that \be\la{bc19}\begin{cases}\sqrt{\n^{R_j}_0}w^{R_j}_0 \rightharpoonup \sqrt{\n_0} w_0 \mbox{ weakly in } L^2(\rr) ,\\
\na w_0^{R_j}\rightharpoonup \na w_0 \mbox{ weakly in } H^1(\rr).\end{cases}\ee It follows from    \eqref{bb12}   that $w_{0}^{R}$ satisfies
 \be \la{bc12} \ba &-\mu \lap w_{0}^{R}-\na \left((\mu+\lambda(\n_{0}^{R}))\div w_{0}^{R}\right)+\na( P(\n_{0}^{R})\vp_R)+ R^{-1} w_0^R\\&  =-\n_0^R w_0^R+\sqrt{\n_{0}^{R}} h^R\vp_R +R^{-1}F^R,\ea \ee
 with $\|F^R\|_{L^2(\rr)}\le C $ due to \eqref{2.9-4} and \eqref{2.9-3}. Thus, one can deduce    from \eqref{bc12},  \eqref{bc10}, and \eqref{bc19} that $w_0  $ satisfies
$$\mu \lap w_{0} -\na \left((\mu+\lambda(\n_0))\div w_{0} \right)+\na P(\n_{0}  ) +\n_0  w_0 =\n_0 u_0+\sqrt{\n_{0} } g,$$
 which combined with \eqref{co2} yields that \be\la{baq1} w_0=u_0.\ee

 Next, we get from \eqref{bc12} and \eqref{co2} that \bnn\ba & \limsup\limits_{R_j\rightarrow \infty}\int_{\rr}\left(\mu|\rot w_0^{R_j}|^2 + (2\mu+\lambda(\n_{0}^{R_j}))(\div w_0^{R_j})^2 + \n_0^{R_j}|w_0^{R_j}|^2\right)dx\\&\le \int_{\rr}\left(\mu|\rot u_0 |^2 + (2\mu+\lm(\n_0))(\div u_0 )^2 + \n_0 |u_0 |^2\right)dx,\ea\enn which together with \eqref{bc19}  implies \bnn \lim\limits_{R_j\rightarrow \infty}\int_{\rr}|\na w_0^{R_j} |^2dx=\int_{\rr}|\na u_0 |^2dx,\,\,\lim\limits_{R_j\rightarrow \infty}\int_{\rr} \n_0^{R_j}  |w_0^{R_j} |^2dx=\int_{\rr}\n_0 | u_0 |^2dx.\enn
 This, along with  \eqref{baq1} and \eqref{bc19},  yields that \be \la{3.74} \lim\limits_{R\rightarrow \infty}\left(\|\na ( w_0^{R }-u_0)\|_{L^2(\rr)}+\|\sqrt{\n_0^R}  w_0^{R }-\sqrt{\n_0}u_0 \|_{L^2(\rr)}\right)=0.\ee
 Similar to \eqref{3.74}, we can obtain that \bnn \lim\limits_{R\rightarrow \infty} \|\na^2 ( w_0^{R }-u_0)\|_{L^2(\rr)} =0. \enn

Finally, in terms of Lemma \ref{th0}, the initial-boundary-value problem  \eqref{b2} with the initial data $(\n_0^R,u_0^R)$ has  a classical solution   $(\n^{R},u^{R}) $ on $B_{R}\times [0,T_R].$ For $\ti g$ defined by \eqref{gq1} with $(\n_0,u_0)$ being replaced by $(\n_0^R,u_0^R)$, it follows from \eqref{bb12},  \eqref{2.9-4}, and \eqref{co2} that \bnn\|\ti g\|_{L^2(\O)}\le C, \enn for some  $C$ independent of $R.$ Hence,   there exists a generic positive constant $C$ independent of $R$ such that all those estimates stated in Proposition \ref{pro} and  Lemmas \ref{lem4.5}-\ref{4.m44} hold for $(\n^R,u^R).$
Extending   $( \n^{R},u^{R})$  by zero on $\rr  \setminus B_{R}$  and denoting  $$\nr\triangleq (\vp_R)^{4/\ti\beta}\n^R,\quad w^R\triangleq \vp_Ru^R,$$ with $\vp_R$ as in \eqref{vp1} and $\ti\beta=\min\{\beta,1\},$ we deduce from  \eqref{o1}     and   Lemmas \ref{lem4.5}-\ref{4.m44} that the sequence
$(\ti\n^R,w^R)$ converges weakly, up to the extraction of subsequences, to some limit $(\n,u)$   satisfying  \eqref{1.10}, \eqref{l1.2}, and \eqref{1.a10}. Moreover, standard arguments yield that $(\n,u)$   in fact is a strong solution to   the problem \eqref{n1}-\eqref{n4}. The proof of Theorem  \ref{t2} is completed.

\end{document}